\documentclass[10pt]{amsart}
\usepackage{amssymb,amsmath,amscd,amsfonts,amsthm,mathrsfs}
\usepackage{verbatim, enumerate} 
\usepackage{graphicx, tikz, ifthen}
\usetikzlibrary{calc}
\tikzset{math3d/.style={x={(1cm,0cm)},z={(0.2cm,1.3cm)},y={(0cm,1cm)}}}
\def\seedrawings{1}%mettre la valeur 1 pour les voir et 0 sinon 
\input xy
\xyoption{all}
\usepackage[all]{xy}  
\usepackage{color}

\usepackage[latin1]{inputenc}
\newcommand{\lint}{[\hspace*{-0.06cm}[}
\newcommand{\rint}{]\hspace*{-0.06cm}]}

\newcommand{\N}{\mathbb{N}}
\newcommand{\Z}{\mathbb{Z}}

\newcommand{\R}{\mathbb{R}}
\newcommand{\C}{\mathbb{C}}

\newcommand{\Bc}{\mathcal{B}}
\newcommand{\Cc}{\mathcal{C}}
\newcommand{\Dc}{\mathcal{D}}
\newcommand{\Ec}{\mathcal{E}}
\newcommand{\Gc}{\mathcal{G}}
\newcommand{\Hc}{\mathcal{H}}

\newcommand{\Jc}{\mathcal{J}}
\newcommand{\Kc}{\mathcal{K}}
\newcommand{\Mcc}{\mathcal{M}}
\newcommand{\Lc}{\mathcal{L}}

\newcommand{\Vc}{\mathcal{V}}
\newcommand{\supp}{\mathrm{supp}}  

\newcommand{\Nr}{\mathscr{N}}

\newcommand{\Qr}{\mathscr{Q}}
\newcommand{\Vr}{\mathscr{V}}
\def\buildh#1^#2{\mathrel{\mathop{\kern 0pt#1}\limits^{#2}}}
\def\build#1_#2^#3{\mathrel{\mathop{\kern 0pt#1}\limits_{#2}^{#3}}}
\newtheorem{theorem}{Theorem}[section]%{\bf}{\it }
\newtheorem{proposition}{Proposition}[section]%{\bf}{\it }
\newtheorem{lemma}{Lemma}[section]%{\bf}{\it }
\newtheorem{definition}{Definition}[section]%{\bf}{\it }
\newcommand{\wlim}{{\rm w}-\lim}

\numberwithin{equation}{section}

\def \bone{\mathbf{1}}

\def \rmi{{\rm i}}

\title{The essential spectrum of the discrete Laplacian on Klaus-sparse graphs}
\begin{document}
%%%%%%%%%%%%%%%%%%%%%%%AUTHOR%%%%%%%%%%%%%%%%%%%%%%%%%%%%%%%%%%%%%%%%%%%
\author{Sylvain Gol\'enia}
\address{Univ. Bordeaux, CNRS, Bordeaux INP, IMB, UMR 5251,  F-33400, Talence, France}
\email{sylvain.golenia@math.u-bordeaux.fr}
\author{Fran\c coise Truc}
\address{Institut Fourier, UMR 5582 du CNRS Universit\'e de Grenoble I, BP 74, 38402 Saint-Martin d'H\`eres, France}
\email{francoise.truc@univ-grenoble-alpes.fr}
\date{Version of \today}
%Abstract%%%%%%%%%%%%%%%%%%%%%%%%%%%%%%%%%%%%%%%%%%%%%%%%%%%%%%%%%%%%%%%
\begin{abstract}
In 1983, Klaus studied a class of potentials with bumps and computed the essential spectrum of the associated Schr\"odinger operator with the help of some localisations at infinity. A key hypothesis is that the distance between two consecutive bumps tends to infinity at infinity. In this article, we introduce a new class of graphs (with patterns) that mimics this  situation, in the sense that the distance between two  patterns tends to infinity at infinity. These patterns tend, in some way, to asymptotic graphs. They are the localisations at infinity. Our result is that the essential spectrum of the Laplacian acting on our graph is given by the union of the spectra of the Laplacian acting on 
the asymptotic graphs. We also discuss the question of the stability of the essential spectrum in the appendix.
\end{abstract}

\maketitle

\tableofcontents
\section{Introduction}
The computation of the essential spectrum of an operator is a standard question in spectral theory. For a large family of Schr\"odinger operators, it is well-known that the essential spectrum is characterised by the behaviour at infinity of the potential. In 1983, Klaus introduces in his article \cite{Kla} a type of potential with bumps with a crucial feature, that is, the distance between two such bumps tends to infinity. He computes the essential spectrum of this $1d$ (continuous) Schr\"odinger operator in terms of the union of the spectrum of some simpler operators. The common patterns, that define the \emph{localisations at infinity},  are given by the behaviour of the potential. In \cite{LaSi}, this notion is illustrated by  the use of \emph{R-limits}, see also \cite{CWL} for this concept and references therein. The example of Klaus was generalised and encoded in some $C^*$-algebraic context in \cite{GeIf1, GeIf2}, see also \cite{MaPuRi}. We refer to \cite{Ge} for more general results and historical references. We mention also \cite{NaTa1, NaTa2} for recent developments in a sparse context. 

In the context of graphs, the computation of essential spectra is done in many places e.g., \cite{Kel, Ra, SaSu}. In \cite{BrDeEl, El} they extend the $R$-limit technique to discrete graphs. We refer to \cite{Gol1, GeGo1} for a $C^*$-algebra approach.

Our motivation in this paper is to analyse a graph analog of the example of Klaus where the ``bumps" are no more due to a potential but to patterns coming from the structure of the graph. We call this family of graphs \emph{Klaus-sparse graphs}. The approaches of $R$-limits and $C^*$-algebras do not seem to apply here. We rely directly on the construction of Weyl sequences.

We start with some definitions and fix our notation for graphs. We refer to
\cite{CdV, Chu} for surveys on the matter.
Let $\Vc$ be a countable set. Let $\Ec:=\Vc\times \Vc\rightarrow [0,\infty)$
and assume that 

\[\Ec(x,y)=\Ec(y,x), \quad \mbox{ for all } x,y\in \Vc.\]

We say that $\Gc:=(\Vc,\Ec, m)$ is an
unoriented weighted graph with \emph{vertices}~$\Vc$ and
\emph{weighted edges}~$\Ec$ and where $m$ is a positive weight on the vertices
\[m:\Vc\rightarrow
(0,\infty).\]
In the setting of electrical networks,
the weights correspond to the conductances.
We say that $x,y\in \Vc$ are \emph{neighbours} if $\Ec(x,y)\neq 0$
and denote it by $x\sim y$. We say that there is a \emph{loop} in $x\in \Vc$
if $\Ec(x,x)\neq 0$.  The set of \emph{neighbours} of $x\in \Ec$
is denoted by 
\[\Nr_\Gc(x):=\{y\in \Ec,  x\sim y\}.\]
A graph is \emph{locally finite} if $|\Nr_\Gc(x)|$ is finite for all $x\in
\Vc$. 
%Finally, as we are dealing with magnetic fields, we fix a
%phase 
%
%\[\theta:\Vr\times \Vr\rightarrow [-\pi, \pi], \mbox{ such that }
%\theta(x,y)= - \theta(y,x).\]
%
%We set $\theta_{x,y}:=\theta(x,y)$.
   A graph is \emph{connected}, 
if for all $x,y\in \Vc$, there exists an $x$-$y$-\emph{path}, i.e.,
there is a finite sequence
\[(x_1,\dotsc,x_{N+1})\in \Vc^{N+1} \mbox{ such that }
x_1=x, \, x_{N+1}=y \mbox{ and } x_n\sim x_{n+1}, 
\]
%\[ (x_1,\dotsc,x_{N+1})\in \Vc^{N+1} \mbox{ such that } x_1=x, \, x_{N+1}=y \mbox{ and } x_n\sim x_{n+1}, \]
%for all $n\in\{1,\dotsc,N\}$. 
%The minimal possible $N$ is called the
%(unweighted) \emph{distance} between  $x$ and $y$.

We recall that a graph~$\Gc$ is \emph{simple} if
$\Ec$ has values in $\{0,1\}$, $m=1$,  and if it has no loop. 
%A \emph{bi-partite} graph is a graph whose vertex set can be partitioned
%into two subsets in such a way that no two points in the same subset
%are neighbors. Trees are bi-partite graphs.
\vspace*{0.5cm}
\begin{center}
{\bf In the sequel, all graphs $\Gc=(\Vc, \Ec, m)$ are  locally
finite, connected and have no loop.}
\end{center}
\vspace*{0.5cm}

We now associate a certain Hilbert space and some operators on it to a given 
graph $\Gc=(\Vc, \Ec, m)$. Let $\ell^2_m(\Vc):=\ell^2(\Vc, m;
\C)$ be the set of functions $f:\Vc\rightarrow \C$, such that
$\|f\|^2_{\ell^2_m(\Vc)}:=\sum_{x\in   \Vr} m(x)|f(x)|^2$ is finite. The associated
scalar product is given by 
\[\langle f, g\rangle_{\Gc}:= \langle f, g\rangle_{\ell^2_m(\Vc)}:= \sum_{x\in \Vc}  
m(x) \overline{f(x)}g(x),  \text{ for } f,g\in \ell^2_m(\Vc).\]
We also denote by $\Cc_c(\Vc)$ the set of functions $f:\Vc\rightarrow
\C$, which have finite support. We define the quadratic form:
\begin{align}\label{e:quadra}
\Qr_\Gc(f,f):=\frac{1}{2}\sum_{x,y\in \Vc}
\Ec(x,y)|f(x)-f(y)|^2\geq 0, \mbox{ for 
} f\in \Cc_c(\Vc).
\end{align}
It is closable and there exists a unique self-adjoint operator
$\Delta_\Gc$, such that  
\[\Qr_\Gc(f,f)= \langle f,
\Delta_\Gc f\rangle_\Gc, \mbox{ for 
} f\in \Cc_c(\Vc)
\]
and $\Dc(\Delta_{\Gc}^{1/2})=
\Dc(\Qr_{\Gc})$, where the latter is the completion of
$\Cc_c(\Vc)$ under $\|\cdot\|^2 + \Qr_{\Gc}(\cdot,
\cdot)$. This operator is the \emph{Friedrichs extension} associated
to the form $\Qr_{\Gc}$ (e.g., \cite{Gol3}). It acts as follows: 
\begin{align}\label{e:Deltam}
\Delta_{\Gc} f(x) :=\displaystyle \frac{1}{m(x)}\sum_{y\in \Vc} \Ec(x,y)
(f(x)- f(y)), \mbox{ for
} f\in \Cc_c(\Vc).
\end{align}
When $m=1$ and $\Ec$ has values in $\{0,1\}$, the operator is 
essentially self-adjoint on $\Cc_c(\Vc)$, c.f., \cite{Woj}. A large literature is 
devoted to this subject. 

We define the \emph{degree} associated to $\Gc= (\Vc, \Ec, m)$ by
\[\deg_\Gc(x):=\frac{1}{m(x)} \sum_{y\in \Vc} 
\Ec(x,y), \text{ for } x\in \Vc.\] 
Given a function $V:\Vc\to \C$, we
denote  by $V(\cdot)$ the operator of multiplication by $V$. It is
elementary that $\Dc(\deg_\Gc^{1/2}(.))\subset\Dc(\Delta^{1/2}_{\Gc
  })$. Indeed, one has: 
\begin{align}\nonumber
0\leq \langle f, \Delta_{\Gc} f\rangle_{\Gc}&= \frac{1}{2} \sum_{x\in \Vc} \sum_{y\sim
  x} \Ec(x,y)|f(x)-f(y)|^2
\\
\label{e:majo}
&\leq \sum_{x\in \Vc} \sum_{y\sim 
  x} \Ec(x,y)(|f(x)|^2+|f(y)|^2) = 2\langle f, \deg_\Gc(\cdot) f\rangle_{\Gc},
\end{align}
for $f\in\Cc_c(\Vc)$. Moreover,  setting $\tilde \delta_x(y):=  m^{-1/2}(x)\bone_{\{x\}}(y)$ for any $x,y \in \Vc$, $\langle \tilde \delta_x, \Delta_{\Gc} \tilde \delta_x\rangle = \deg_\Gc(x)$,
so $\Delta_{\Gc}$ is bounded if and only if $\sup_{x\in
  \Vc}\deg_\Gc(x)$ is finite, e.g.\ \cite{KL, Gol3}. Here we have used the following standard notation: for any given set $X$, $\bone_X(x):=1$ is $x\in X$ and $0$ otherwise.

The aim of our work is to study the essential spectrum of $\Delta_{\Gc}$ for  Klaus-sparse graphs. The precise definition is given in Section \ref{s:KS}
but let us give a rough description of such a graph: it  consists of a (double infinite) family of finite graphs  $\{\Gc_{i,k}\}_{(i,k)}$ which are patterns that are connected by a ``medium graph" $\Gc_\Mcc$. The latter has a uniformly bounded degree. For any given $k$, there is a kind of increasing limit of $\Gc_{i,k}$ when $i\to \infty$ which is an infinite graph $\Gc_{\infty,k}$. The graphs 
$\{\Gc_{\infty,k}\}_{k}$ are the localisations at infinity of $\Gc$. Moreover the distance in $\Gc$ between the  $\Gc_{i,k}$ goes to infinity as $\lim (i,k)\to \infty$ with respect to  the Fr\'echet filter. We prove the following theorem~:

\begin{theorem}\label{t:simpleF}
Assume that $\Gc:= (\Vc, \Ec, m)$ is  a Klaus-sparse graph. Using notation of Definition \ref{d:KS2}, we have
\begin{enumerate}[1)]
\item $\Delta_\Gc$ is essentially self-adjoint on $\Cc_c(\Vc)$.
\item The essential spectrum of $\Delta_\Gc$ is given by the union of the spectra of $\Delta_{\Gc_{\infty, k}}$,  the localisations at infinity, namely:
\[\sigma_{\rm ess} (\Delta_\Gc)= \overline{ \bigcup_{k\in \Jc}\sigma(\Delta_{\Gc_{\infty, k}})}.\]
\end{enumerate}
\end{theorem}
First of all, to simplify the presentation we do not include any potential and stick to perturbation of graphs. The proof of 1) is inspired by \cite{Gol2}. Let us mention that the inclusion $\subset$ in 2) holds for a wider class of graphs and is the easy part of 2). The reverse inclusion is the interesting and difficult part. 
To our knowledge, in the context of $C^*$-algebra, the only result that could tackle this  issue is in \cite{Ge}. However it is complicated to compare his results to ours, since  they are given in terms of some abstract ultra-filters. Moreover, in our result the union that we obtain is minimal, e.g., Section \ref{s:sharpunion}.  We now compare to \cite{El}. In his chapter 4, the author's more general  results  overlap with ours but they do not contain ours. They are complementary.  Let us explain in which sense they do not contain our result. First we point out that the author is dealing only with bounded operators and with $m=1$. We do not believe that it is a strong obstacle for his method. %However, he is assuming that there is $\alpha, \alpha'>0$, and $c, c'$ such that
%\[c r^\alpha \leq |B_\Gc(x, r)|\leq c' r^{\alpha'}, \quad \forall x\in \Vc.\]
However, his approach relies fundamentally on a reverse Shnol's Theorem that is fulfilled under a uniform sub-exponential growth: $\forall \gamma>1, \exists C>0$ so that $\forall r\in\N^*$
\begin{align}\label{e:growth}
\sup_{x\in \Vc}| \{y\in \Vc, d_\Gc(x,y)=r\}|\leq C \gamma^r.
\end{align}
We construct in Section \ref{s:nonex} a general family of Klaus-sparse graphs that do not have a uniform exponential growth.

Finally in the appendix, Theorem \ref{t:stabess} ensures the stability of the essential spectrum under a perturbation of the metric that is small at infinity. 

{\bf Acknowledgements: } We would like to thank Eric Amar and Michel Bonnefont for fruitful discussions. We thank Latif Eliaz and Jonathan Breuer for mentioning their interesting results.

\section{Klaus-sparse graphs}\label{s:KS}

\subsection{Further notation}\label{s:FN} We introduce some further  notation and definitions. Given $a,b\in \Z$, we denote by $\lint
a,b\rint:=[a,b]\cap \Z$ and  $\lint
a,\infty\lint:=[a,\infty[\cap \Z$. Given $X\subset Y$, we denote  by $X^c:= Y\setminus X$ the complementary set of $X$, when no confusion can arise. Given $\Hc, \Kc$ Hilbert spaces, we denote by $\Bc(\Hc, \Kc)$ and $\Kc(\Hc, \Kc)$ the set of bounded and compact operators from $\Hc$ to $\Kc$. Set also $\Bc(\Hc):=\Bc(\Hc, \Hc)$ and $\Kc(\Hc):=\Kc(\Hc, \Hc)$.

Given $\Gc_i:= (\Vc_i,\Ec_i, m_i)$ and $x_i\in \Vc_i$  for $i=1,2$, we say that
$\Gc_1$ is \emph{induced} by $\Gc_2$ and denote it by $\Gc_1\subset \Gc_2$ if
there is an injection $f:\Vc_1\to \Vc_2$ such that 
\[f(\Vc_1)\subset \Vc_2, \quad \Ec_1(x, y) =
\Ec_2 (f(x), f(y)), \text{ and } m_1 (x)= m_2(f(x)), \quad \forall x,y\in \Vc_1.\]
We shall  write $(\Gc_1, x_1)\subset (\Gc_2, x_2)$ if  we have in addition $f(x_1)=x_2$. To simplify, we shall  often simply write
\[\Vc_1\subset \Vc_2, \quad \Ec_1=
\Ec_2 |_{\Vc_1, \Vc_1}, \quad  m_1= m_2|_{\Vc_1}, \text{ and } x_1=x_2. \]
Moreover, given a graph $\Gc=  (\Ec, \Vc, m)$ and $X\subset \Vc$, we denote by $[X]^\Gc:= (X, \Ec|_{X\times X}, m  |_{X})$ the induced graph of $\Gc$ by $X$.

We denote by $d_\Gc$ the (unweighted) distance for $\Gc$ (over $\Vc$) given by 
\[d_\Gc(x,y) := \min(n, x=x_0\sim x_1 \sim \ldots \sim x_n=y, \text{ with } x_i\in \Vc),\]
 for $x,y\in \Vc$ when $x\neq y$ and $d_\Gc(x,x):=0$. It is a distance on $\Vc$ when $\Gc$ is connected. Given $r\geq 0$, we set
 \[B_\Gc(x, r):= \{y\in \Vc, d_{\Gc}(x,y)\leq r\}.\]

\subsection{The Klaus-sparse graph} In the introduction we explained by hand-waving the notion of Klaus-sparse graph that we present in this article. Here is the precise definition.

\begin{definition}\label{d:KS2}
A graph $\Gc:=( \Vc,\Ec, m)$ is called a \emph{Klaus-sparse} graph if the
following is satisfied.
\begin{itemize}
\item \underline{Set of patterns}: Let $\Jc$ be at most countable such that $|\Jc|\geq 1$.  There exist a family of connected subgraphs $\Gc_{i,k}:=( \Vc_{i,k},\Ec_{i,k}, m_{i,k})$ with $i\in \N$ and $k\in \Jc$, that are induced by $\Gc$ . 
\item \underline{Localisations at infinity}: There are $\Gc_{\infty,k}:=( \Vc_{\infty,k},\Ec_{\infty,k}, m_{\infty,k})$, for all $k\in \Jc$, and  $\Gc_\Mcc:=(\Vc_\Mcc, \Ec_\Mcc, m_\Mcc)$ and $x_\Mcc\in \Vc_\Mcc$, with uniformly bounded degree:
\begin{align}\label{e:Msup}
\sup_{y\in\Vc_\Mcc} \deg_{\Gc_\Mcc}(y)<\infty.
\end{align}
Here $\Mcc$ stands for \emph{medium}. 
\end{itemize}
These two families verify the following compatibilities:
\begin{enumerate}
\item For all $i\in \N$ and $k\in\Jc$, there are $0< r_{i,k}^{\rm int} <r_{i,k}^{\rm ext}$ such that
\[\lim_{(i,k)\to \infty}r_{i,k}^{\rm int} = \infty \quad \text{ and } \quad \lim_{(i,k)\to \infty} r^{\rm ext}_{i,k} -r^{\rm int}_{i,k} =\infty,\]
where the limit is taken with respect to the Fr\'echet filter, i.e., the one given by the complementary of finite sets. 
\item For all $i\in\N$  and $k\in \Jc$, there exist $x_{i,k}
\in \Vc_{i,k}$  and  $x_{\infty,k}\in \Vc_{\infty,k}$, so that
\[([B_\Gc(x_{i,k}, r_{i,k}^{\rm ext})]^\Gc, x_{i,k})\subset (\Gc_{\infty,k}, x_{\infty, k}),\]
and for all $r>0$ and $k\in \Jc$, there is $i\in \N$ such that
\[([B_{\Gc_{\infty,k}}(x_{\infty, k}, r)]^{\Gc_{\infty,k}}, x_{\infty, k} ) \subset ([B_\Gc(x_{i,k}, r_{i,k}^{\rm int}-1)]^\Gc, x_{i,k}).\] 
\item We set
\[\{\Cc_l\}_{l\in \Lc}:= \{\text{connected component of }[\Vc \setminus \cup_{i\in \N, k \in \Jc} B_\Gc(x_{i,k,} r_{i,k}^{\rm int}-1)]^{\Gc}\}.\]
Here $\Lc\subset \N^*$. It is finite or not. We suppose that for all $l\in \Lc$, there are a vertex $x_l$ of $\Cc_l$ and an order $\lhd$ on $\Lc$ such that 
\[\forall l,m\in \Lc,\, l\lhd  m \implies (\Cc_l, x_l) \subset (\Cc_m, x_m)\subset (\Gc_\Mcc, x_\Mcc).\]

%and such that the increasing limit of $\{\Cc_i\}_{i\in I}$, with respect to $\leq$, is equal to $\Gc_\Mcc$. \blue{needs to be precised}
%%%
%\item[($d_{\rm strong}$)]  \red{changement, hypoth\`ese plus forte}
%\[\lim_{i\to \infty} \inf_{j\in \N\setminus\{i\}}d_\Gc(B_\Gc(x_i, r_i^{\rm ext}), B_\Gc(x_j, r_j^{\rm ext})) =\infty\]
%\item[($d$)] For all $i\in \N$, \red{ancienne hypoth\`ese }
%\[\lim_{j\to \infty} d_\Gc(B_\Gc(x_i, r_i^{\rm ext}), B_\Gc(x_j, r_j^{\rm ext})) =\infty\]
\item[($d$)]  For all $i,j\in \N$ and $k,k'\in \Jc$, such that $(i,k)\neq (j, k')$, we have 
\[B_\Gc(x_{i,k}, r_{i,k}^{\rm ext}) \cap B_\Gc(x_{j,k'}, r_{j,k'}^{\rm ext})=\emptyset.\]
\item[($e$)] For all $r>0$, there is $(i,k)\in \N\times \Jc$ such that
\begin{align}\label{e:backside}
[B_{\Gc_\Mcc}(x_\Mcc, r)]^{\Gc_\Mcc} \subset  [B_\Gc(x_{i,k}, r_{i,k}^{\rm ext}-1)\setminus  B_\Gc(x_{i,k}, r_{i,k}^{\rm int}+1)]^{\Gc}.
\end{align}%and
%\[B_\Gc(x_{i,k}, r_{i,k'}^{\rm ext}) \cap B_\Gc(x_{i,k}, r_{i,k'}^{\rm ext})=\emptyset.\]
\end{enumerate}
\end{definition}

At first sight, it looks a bit abstruse but will be more intuitive by relying on the figures. In Figure \ref{fig_antitree}, $C_i\simeq \lint -a_i, b_i \rint$ where $a_i$ and $b_i$ are integers that tend to $\infty$. The medium graph is $\Gc_\Mcc = \Z$, see Figure \ref{fig_antitree_inf}. Note that we could have chosen $\Gc_\Mcc=\N$.
 
In Figure \ref{fig_Z2like},  $[\Vc \setminus \cup_{i\in \N} B_\Gc(x_i, r_i^{\rm int})]^{\Gc}$ is connected. Therefore $\Lc$ is reduced to an element and $\Gc_\Mcc= \Z^2$, see Figure \ref{fig_Z2like2}.

%%%%%%%%%%%%%%%%%
\if\seedrawings1
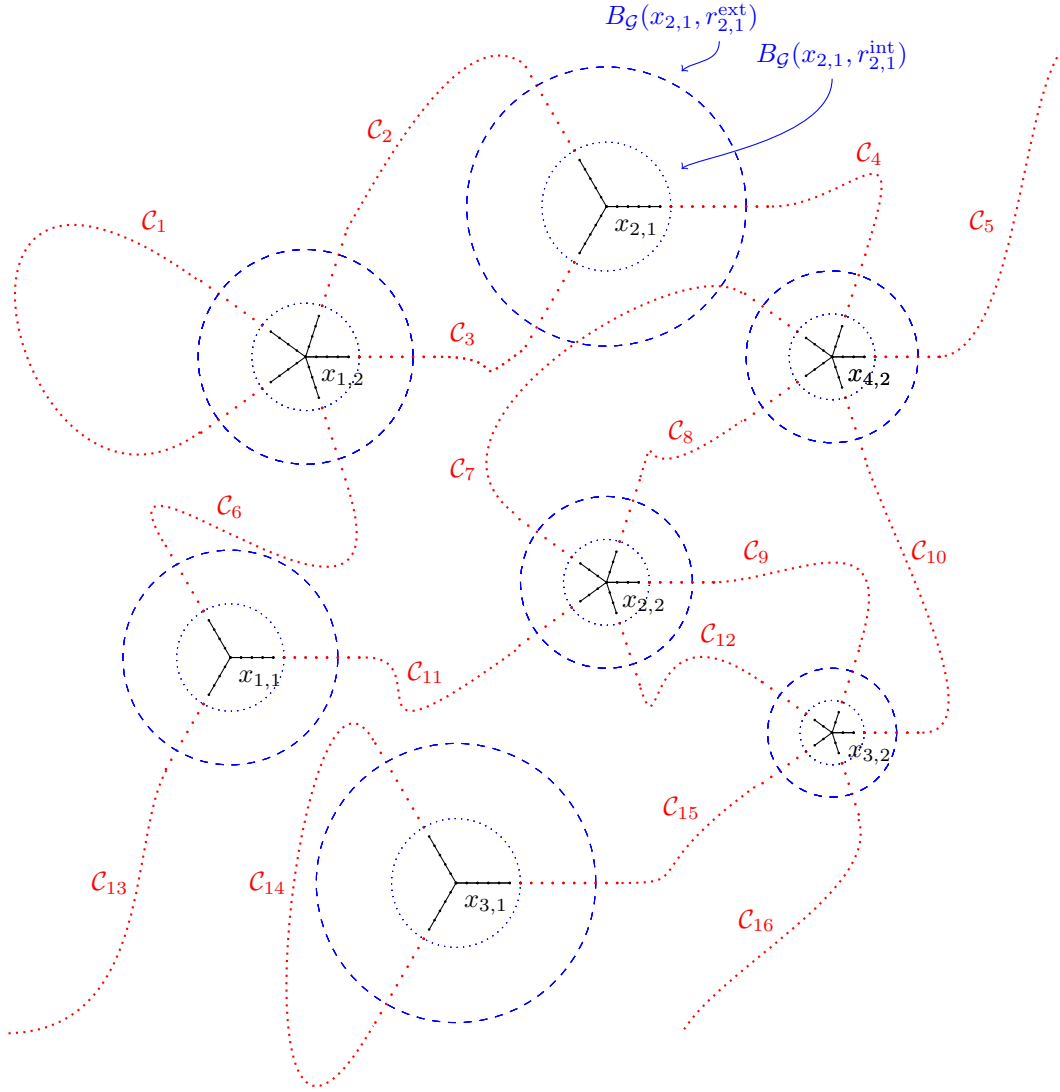
\begin{figure}
  \begin{tikzpicture}
  %\draw [color=green!20,very thin, step=0.2] (-5,-5) grid (10.2,10.2);
  %\draw [color=green,very thin, step=1] (-5,-5) grid (10.2,10.2);
  %\filldraw (0,0) circle(0.1) [color =green];
  \def\xab{  7/0/3/8/5, 7/5/4/10/5, 4/2/4/10/5, -1/1/5/12/3,  0/5/5/12/5, 4/7/6/15/3, 2/-2/6/15/3} % x / y/ rint/ nb de point / nb de branches
  \def\scal{7}
  \foreach \x/\y/\rint/\r/\Ang in \xab {
    \foreach \ang in {0,...,\Ang}{
    \draw[-] (\x,\y) -- +(\ang*360/\Ang:{(\rint-1)/\scal}) ;
    \draw [dotted, color=blue] (\x,\y) circle({\rint/\scal}) ;
    \draw [dashed, color=blue] (\x,\y) circle({(\r-2)/\scal}) ;
      \foreach \s in {0,...,\rint} { 
       \filldraw (\x,\y) + (\ang*360/\Ang:\s/\scal) circle (0.1mm) [color = black] ;
      } 
            \foreach \s in {\rint,...,\r} { 
       \filldraw (\x,\y) + (\ang*360/\Ang:\s/\scal) circle (0.1mm) [color = red] ;
      } 
    }
   }
   \draw[dotted, thick, color=red] (0.55,6.7) .. controls +(60:1) and +(144:1) .. (2.9,8.9);
    \draw[dotted, thick, color=red] (1.8,5) .. controls +(0:1) and +(-144:1) .. (2.9,5.1);
     \draw[dotted, thick, color=red] (-1.4,4) .. controls +(-144:3) and +(144:5) .. (-1.4,6);
      \draw[dotted, thick, color=red] (-1.9,2.55) .. controls +(120:2) and +(-72:3) .. (0.55,3.30);
       \draw[dotted, thick, color=red] (1.8,5) .. controls +(0:1) and +(-144:1) .. (2.9,5.1);
        \draw[dotted, thick, color=red] (.8,1) .. controls +(0:1) and +(-144:3) .. (2.8,1.1);
        \draw[dotted, thick, color=red] (5.8,5.85) .. controls +(144:1) and +(144:2) .. (2.8,2.9);
        \draw[dotted, thick, color=red] (-1.85,-.5) .. controls +(-120:1) and +(0:2) .. (-4,-4);
        \draw[dotted, thick, color=red] (0.9,-.1) .. controls +(120:2) and +(-120:4) .. (0.9,-3.9);
        \draw[dotted, thick, color=red] (4.2,-2) .. controls +(0:1) and +(-144:2) .. (6,-.7);
        \draw[dotted, thick, color=red] (4.45,.6) .. controls +(-72:1) and +(144:2) .. (6.1,0.65);
        \draw[dotted, thick, color=red] (4.45,3.4) .. controls +(72:1) and +(-144:2) .. (5.8,4.15);
        \draw[dotted, thick, color=red] (5.4,2) .. controls +(0:1) and +(72:2) .. (7.35,1.1);
        \draw[dotted, thick, color=red] (6.2,7) .. controls +(0:1) and +(72:2) .. (7.45,6.4);
       \draw[dotted, thick, color=red] (7.35,-1.15) .. controls +(-72:1) and +(60:1) .. (5,-4);
       \draw[dotted, thick, color=red] (8.4,5) .. controls +(0:1) and +(-120:1) .. (10,9);
       \draw[dotted, thick, color=red] (8.2,0) .. controls +(0:1) and +(-72:1) .. (7.45,3.6);
       \node at (-0.6,0.7) {$x_{1,1}$};
       \node at (4+0.4,7-0.3) {$x_{2,1}$};
       \node at (2+0.4,-2-0.3) {$x_{3,1}$};
       \node at (0+0.5,5-0.3) {$x_{1,2}$};
       \node at (4+0.5,2-0.3) {$x_{2,2}$};
       \node at (7+0.5,0-0.3) {$x_{3,2}$};
       \node at (7+0.5,5-0.3) {$x_{4,2}$};
       \node at (7+0.5,5-0.3) {$x_{4,2}$};
       \node[color=red] at (-2,6.8) {$\Cc_1$};
       \node[color=red] at (1,8) {$\Cc_2$};
       \node[color=red] at (2.1,5.3) {$\Cc_3$};
       \node[color=red] at (7.5,7.7) {$\Cc_4$};
       \node[color=red] at (9,6.8) {$\Cc_5$};
       \node[color=red] at (-1,3) {$\Cc_6$};
       \node[color=red] at (2.1,3.5) {$\Cc_7$};
       \node[color=red] at (5,4) {$\Cc_8$};
       \node[color=red] at (6,2.4) {$\Cc_9$};
       \node[color=red] at (8.3,2.4) {$\Cc_{10}$};
       \node[color=red] at (1.6,0.8) {$\Cc_{11}$};
       \node[color=red] at (5.5,1.3) {$\Cc_{12}$};
       \node[color=red] at (-2.6,-2) {$\Cc_{13}$};
       \node[color=red] at (-0.5,-2) {$\Cc_{14}$};
       \node[color=red] at (5,-1) {$\Cc_{15}$};
       \node[color=red] at (6,-2.5) {$\Cc_{16}$};
       \node[color=blue] at (5,9.5) {$B_\Gc(x_{2,1}, r^{\rm ext}_{2,1})$};
       \draw[->, color=blue] (5.5,9.2)  .. controls +(-90:.3) and +(45:.2) .. (5,8.8);
       \node[color=blue] at (7,9) {$B_\Gc(x_{2,1}, r^{\rm int}_{2,1})$};
        \draw[->, color=blue] (7,8.7)  .. controls +(-90:.8) and +(45:.2) ..  (5,7.5);
  \end{tikzpicture}
  \caption{$\Gc$ is a Klaus-sparse star-like graph. Here $\Lc=\N$.}
  \label{fig_antitree}
\end{figure}

\begin{figure}
  \begin{tikzpicture}
  %\draw [color=green!20,very thin, step=0.2] (-3,-3) grid (10.2,3.2);
  %\draw [color=green,very thin, step=1] (-3,-3) grid (10.2,3.2);
  %\filldraw (0,0) circle(0.1) [color =green];
  \def\xab{  -2/0/4/16/3, 2.5/0/4/16/5, 7.5/0/4/16/2 } % x / y/ rint/ nb de point / nb de branches
  \def\scal{7}
  \foreach \x/\y/\rint/\r/\Ang in \xab {
   \filldraw (\x,\y) circle (0.5mm) [color = black] ;
    \foreach \ang in {0,...,\Ang}{
    \draw[-] (\x,\y) -- +(\ang*360/\Ang:{(\r-2)/\scal}) ;
    %\draw [dotted, color=blue] (\x,\y) circle({\rint/\scal}) ;
    \draw [dashed, color=blue] (\x,\y) circle({(\r-4)/\scal}) ;
      \foreach \s in {0,...,\r} { 
       \filldraw (\x,\y) + (\ang*360/\Ang:\s/\scal) circle (0.2mm) [color = black] ;
      } 
    }
   }
  \node at (-1.6,.2) {$x_{\infty, 1}$};
  \node at (3.1,.2) {$x_{\infty, 2}$};
   \node at (7.9,.2) {$x_{\Mcc}$};
    \node at (-2,-2.3) {$\Gc_{\infty,1}$};
    \node at (2.5,-2.3) {$\Gc_{\infty,2}$};
    \node at (7.5,-2.3) {$\Gc_{\Mcc}$};
     \node[color=blue] at (-1,2.1) {$B_{\Gc_{\infty,1}}(x_{\infty,1}, r)$};
      \node[color=blue] at (4.6,1.7) {$B_{\Gc_{\infty,2}}(x_{\infty,2}, r)$};
       \node[color=blue] at (7.5,2.1) {$B_{\Gc_{\Mcc}}(x_\Mcc, r)$};
  \end{tikzpicture}
  \caption{Localisations at infinity for the Star-like graph given in Figure \ref{fig_antitree}}
  \label{fig_antitree_inf}
\end{figure}
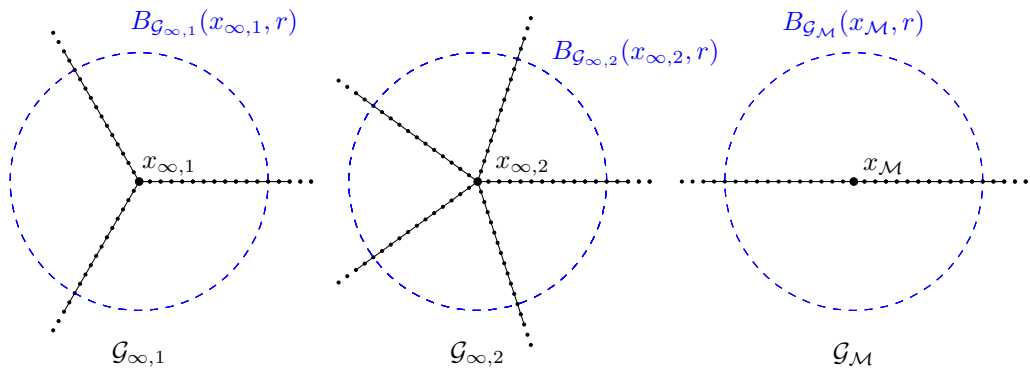

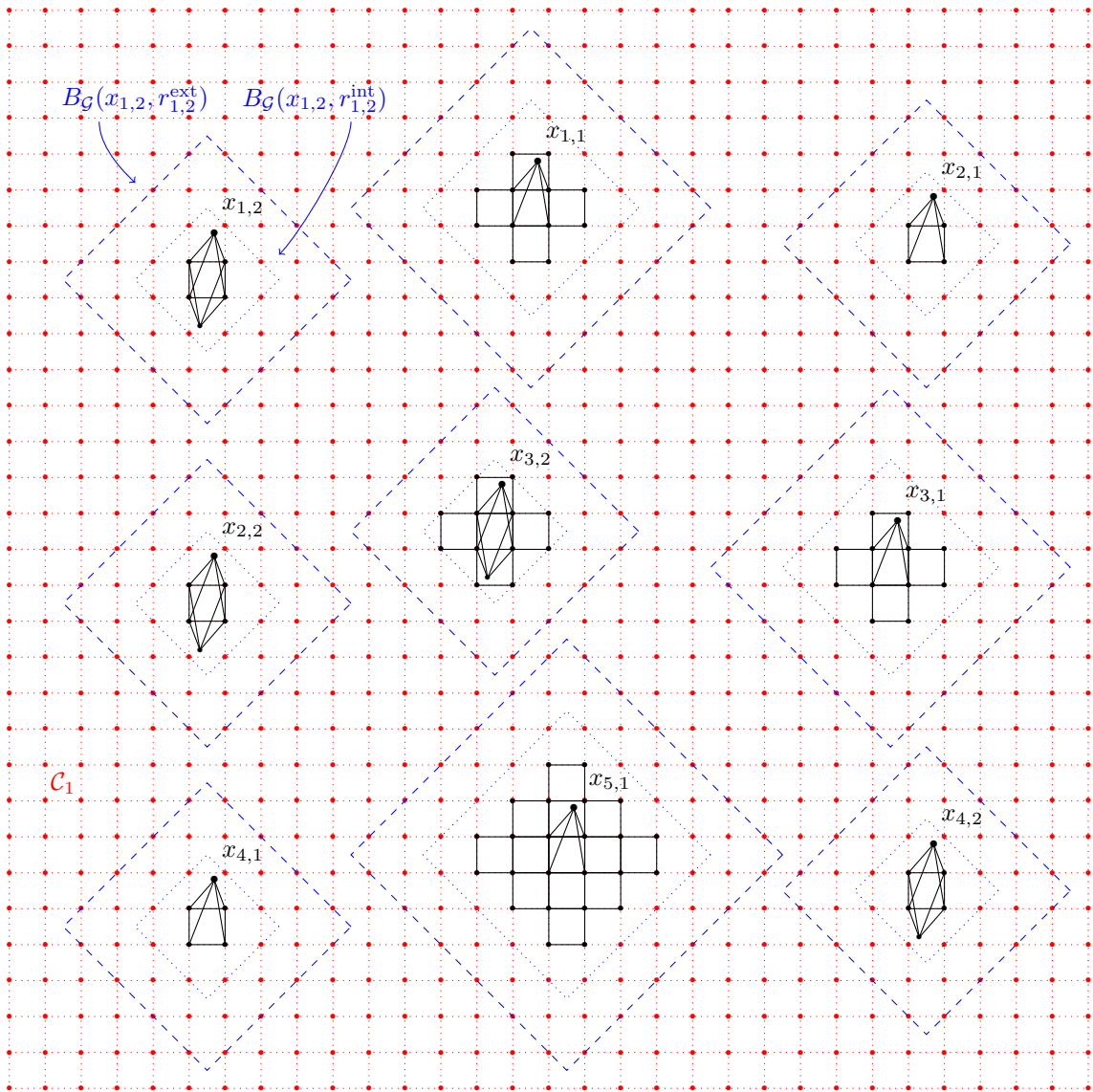
\begin{figure}
\begin{tikzpicture}[math3d, scale=0.5]
	\def \xmax{24}
	%\draw (0,0,0) [color =green] node{$\bullet$};
	 \foreach \x in {-6,...,\xmax}{
	 \draw[color =red!100, dotted] (\x,-6,0) -- (\x, \xmax,0);
	 \draw[color =red!100, dotted] (-6,\x,0) -- (\xmax, \x,0);
	 	\foreach \y in {-6,...,\xmax}{
	\draw (\x,\y,0) [color =red!100] node[scale=.5]{$\bullet$};
		}
         } 
         %pyramidea a un sommet
         \def \xab{-1/-2/2/4, 8/18/3/5, 9/0/4/6, 19/17/2/4, 18/8/3/5}
        \foreach \x/\y/\tmin/\tmax in \xab{
         \draw ({\x+.5},{\y+.5},1) [color =black] node[scale=.7]{$\bullet$};
         \draw[color = blue, dashed] ( {\x+.5+ \tmax}, {\y+.5}, 0) -- ( {\x+.5}, {\y+.5+ \tmax}, 0)-- ( {\x+.5- \tmax}, {\y+.5}, 0) -- ( {\x+.5}, {\y+.5- \tmax}, 0)-- cycle;
         \draw[color = blue, dotted] ( {\x+.5+ \tmax-2}, {\y+.5}, 0) -- ( {\x+.5}, {\y+.5+ \tmax-2}, 0)-- ( {\x+.5- \tmax+2}, {\y+.5}, 0) -- ( {\x+.5}, {\y+.5- \tmax+2}, 0)-- cycle;
         \draw [color=black] ({\x+.5},{\y+.5},1)--({\x},{\y},0);
         \draw [color=black] ({\x+.5},{\y+.5},1)--({\x+1},{\y},0);
         \draw [color=black] ({\x+.5},{\y+.5},1)--({\x},{\y+1},0);
	\draw [color=black] ({\x+.5},{\y+.5},1)--({\x+1},{\y+1},0);
		\foreach \xr in {0, ..., \tmin}{
		\foreach \yr in {0, ..., \tmin}{
		\pgfmathparse{\tmin-\xr-1}\let\a\pgfmathresult
		\ifthenelse{\yr<\a}{
		\draw[color =black!100] ({\x+\xr+1},{\yr+\y+1},0) node[scale=.5]{$\bullet$};
		\draw[color =black!100] ({\x-\xr},{\yr+\y+1},0) node[scale=.5]{$\bullet$};
		\draw[color =black!100] ({\x-\xr},{-\yr+\y},0) node[scale=.5]{$\bullet$};
		\draw[color =black!100] ({\x+\xr+1},{-\yr+\y},0) node[scale=.5]{$\bullet$};
		\draw[color =black!100] ({\x+\xr+1},{-\yr+\y},0)-- ({\x-\xr},{-\yr+\y},0);
		\draw[color =black!100] ({\x+\xr+1},{\yr+\y+1},0)-- ({\x-\xr},{\yr+\y+1},0);
		\draw[color =black!100] ({\x+\xr+1},{-\yr+\y},0)-- ({\x+\xr+1},{\yr+\y+1},0);
		\draw[color =black!100] ({\x-\xr},{-\yr+\y},0) -- ({\x-\xr},{\yr+\y+1},0);
		}{}
}{}
		%
%	}
         }
         }
         %pyramides a deux sommets
                  \def \xab{7/9/3/4, -1/16/2/4, -1/7/2/4, 19/-1/2/4}
        \foreach \x/\y/\tmin/\tmax in \xab{
        %\pgfmathparse{\tmax-\2}\let\tmin\pgfmathresult
         \draw[color = blue, dashed] ( {\x+.5+ \tmax}, {\y+.5}, 0) -- ( {\x+.5}, {\y+.5+ \tmax}, 0)-- ( {\x+.5- \tmax}, {\y+.5}, 0) -- ( {\x+.5}, {\y+.5- \tmax}, 0)-- cycle;
         \draw[color = blue, dotted] ( {\x+.5+ \tmax-2}, {\y+.5}, 0) -- ( {\x+.5}, {\y+.5+ \tmax-2}, 0)-- ( {\x+.5- \tmax+2}, {\y+.5}, 0) -- ( {\x+.5}, {\y+.5- \tmax+2}, 0)-- cycle;
         \draw ({\x+.5},{\y+.5},1) [color =black] node[scale=.7]{$\bullet$};
         \draw ({\x+.5},{\y+.5},-1) [color =black] node[scale=.5]{$\bullet$};
         \draw [color=black] ({\x+.5},{\y+.5},1)--({\x},{\y},0);
         \draw [color=black] ({\x+.5},{\y+.5},1)--({\x+1},{\y},0);
         \draw [color=black] ({\x+.5},{\y+.5},1)--({\x},{\y+1},0);
	\draw [color=black] ({\x+.5},{\y+.5},1)--({\x+1},{\y+1},0);
	\draw [color=black] ({\x+.5},{\y+.5},-1)--({\x},{\y},0);
         \draw [color=black] ({\x+.5},{\y+.5},-1)--({\x+1},{\y},0);
         \draw [color=black] ({\x+.5},{\y+.5},-1)--({\x},{\y+1},0);
	\draw [color=black] ({\x+.5},{\y+.5},-1)--({\x+1},{\y+1},0);
	\foreach \xr in {0, ..., \tmin}{
		\foreach \yr in {0, ..., \tmin}{
		\pgfmathparse{\tmin-\xr-1}\let\a\pgfmathresult
		\ifthenelse{\yr<\a}{
		\draw[color =black!100] ({\x+\xr+1},{\yr+\y+1},0) node[scale=.5]{$\bullet$};
		\draw[color =black!100] ({\x-\xr},{\yr+\y+1},0) node[scale=.5]{$\bullet$};
		\draw[color =black!100] ({\x-\xr},{-\yr+\y},0) node[scale=.5]{$\bullet$};
		\draw[color =black!100] ({\x+\xr+1},{-\yr+\y},0) node[scale=.5]{$\bullet$};
		\draw[color =black!100] ({\x+\xr+1},{-\yr+\y},0)-- ({\x-\xr},{-\yr+\y},0);
		\draw[color =black!100] ({\x+\xr+1},{\yr+\y+1},0)-- ({\x-\xr},{\yr+\y+1},0);
		\draw[color =black!100] ({\x+\xr+1},{-\yr+\y},0)-- ({\x+\xr+1},{\yr+\y+1},0);
		\draw[color =black!100] ({\x-\xr},{-\yr+\y},0) -- ({\x-\xr},{\yr+\y+1},0);
		}{}
}{}
		%
%	}
         }
}
\node at (9.5,20.5) {$x_{1,1}$};
\node at (20.5,19.5) {$x_{2,1}$};
\node at (19.5,10.5) {$x_{3,1}$};
\node at (.5,.5) {$x_{4,1}$};
\node at (10.7,2.5) {$x_{5,1}$};
\node at (.5,18.5) {$x_{1,2}$};
\node at (.5,9.5) {$x_{2,2}$};
\node at (8.5,11.5) {$x_{3,2}$};
\node at (20.5,1.5) {$x_{4,2}$};
\node[color=red] at (-4.5,2.5) {$\Cc_1$};
\node[color=blue] at (-2.5,21.5) {$B_\Gc(x_{1,2}, r^{\rm ext}_{1,2})$};
\node[color=blue] at (2.5,21.5) {$B_\Gc(x_{1,2}, r^{\rm int}_{1,2})$};
\draw[->, color=blue] (-3.5,20.9) .. controls +(-90:.8) and +(135:.2) .. (-2.5,19.2);
\draw[->, color=blue] (3.5,20.9) .. controls +(-90:.8) and +(45:.8) .. (1.5,17.2);
 %\draw[->, color=blue] (7,8.7)  .. controls +(-90:.8) and +(45:.2) ..  (5,7.5);
	\end{tikzpicture}
	  \caption{$\Gc$ is a Klaus-sparse $\Z^2$-like graph. Here $\Lc=\{1\}$.}
  \label{fig_Z2like}
\end{figure}

\begin{figure}
% Debut G_{inf,1}
\begin{tikzpicture}[math3d, scale=0.5]
	%\draw (0,0,0) [color =green] node{$\bullet$};
	 \foreach \x in {-5,...,6}{
	 \draw[color =black!100, dotted] (\x,-5,0) -- (\x, 6,0);
	 \draw[color =black!100, dotted] (-5,\x,0) -- (6, \x,0);
	 }
	 \foreach \x in {-4,...,5}{
	 \draw[color =black!100] (\x,-4,0) -- (\x, 5,0);
	 \draw[color =black!100] (-4,\x,0) -- (5, \x,0);
	 }
	 \foreach \x in {-4,...,5}{
	 	\foreach \y in {-4,...,5}{
	\draw (\x,\y,0) [color =black!100] node[scale=.5]{$\bullet$};
		}
         }   
                 \draw [color=black] (0.5,0.5,1)--(0,0,0);
         \draw [color=black] (.5,.5,1)--(1,0,0);
         \draw [color=black] (.5,.5,1)--(0,1,0);
	\draw [color=black] (.5,.5,1)--(1,1,0);
         \draw[color = blue, dashed] (4.5, 0.5, 0) -- ( 0.5,  4.5, 0)-- ( - 3.5, 0.5, 0) -- ( 0.5,- 3.5, 0)-- cycle;
         \draw (0.5,0.5,1) [color =black!100] node[scale=.7]{$\bullet$};
         \node at (0,-6) {$\Gc_{\infty,1}$};
         \node[scale=.6] at (1.55,1.6) {$x_{\infty,1}$};
%         \node[color=blue] at (-2.5,21.5) {$B_\Gc(x_{1,2}, r^{\rm ext}_{1,2})$};
	\node[color=blue] at (2,6.5) {$B_{\Gc_{\infty,1}}(x_{\infty,1}, r)$};
	\draw[->, color=blue] (3.5,6) .. controls +(-90:.8) and +(45:.2) .. (2.5,3.5);
\end{tikzpicture}
\quad \quad
% Debut G_inf,2
\begin{tikzpicture}[math3d, scale=0.5]
	%\draw (0,0,0) [color =green] node{$\bullet$};
	 \foreach \x in {-5,...,6}{
	 \draw[color =black!100, dotted] (\x,-5,0) -- (\x, 6,0);
	 \draw[color =black!100, dotted] (-5,\x,0) -- (6, \x,0);
	 }
	 \foreach \x in {-4,...,5}{
	 \draw[color =black!100] (\x,-4,0) -- (\x, 5,0);
	 \draw[color =black!100] (-4,\x,0) -- (5, \x,0);
	 }
	 \foreach \x in {-4,...,5}{
	 	\foreach \y in {-4,...,5}{
	\draw (\x,\y,0) [color =black!100] node[scale=.5]{$\bullet$};
		}
         }   
                 \draw [color=black] (0.5,0.5,1)--(0,0,0);
         \draw [color=black] (.5,.5,1)--(1,0,0);
         \draw [color=black] (.5,.5,1)--(0,1,0);
	\draw [color=black] (.5,.5,1)--(1,1,0);
         \draw[color = blue, dashed] (4.5, 0.5, 0) -- ( 0.5,  4.5, 0)-- ( - 3.5, 0.5, 0) -- ( 0.5,- 3.5, 0)-- cycle;
         \draw (0.5,0.5,1) [color =black!100] node[scale=.7]{$\bullet$};
         \node at (0,-6) {$\Gc_{\infty,2}$};
         \node[scale=.6] at (1.55,1.6) {$x_{\infty,2}$};
                          \draw [color=black] (0.5,0.5,-1)--(0,0,0);
         \draw [color=black] (.5,.5,-1)--(1,0,0);
         \draw [color=black] (.5,.5,-1)--(0,1,0);
	\draw [color=black] (.5,.5,-1)--(1,1,0);
	\draw (0.5,0.5,-1) [color =black!100] node[scale=.5]{$\bullet$};
	\node[color=blue] at (2,6.5) {$B_{\Gc_{\infty,2}}(x_{\infty,2}, r)$};
	\draw[->, color=blue] (3.5,6) .. controls +(-90:.8) and +(45:.2) .. (2.5,3.5);
\end{tikzpicture}
%debut G_M
\begin{tikzpicture}[math3d, scale=0.5]
	%\draw (0,0,0) [color =green] node{$\bullet$};
	 \foreach \x in {-6,...,6}{
	 \draw[color =black!100, dotted] (\x,-6,0) -- (\x, 6,0);
	 \draw[color =black!100, dotted] (-6,\x,0) -- (6, \x,0);
	 }
	 \foreach \x in {-5,...,5}{
	 \draw[color =black!100] (\x,-5,0) -- (\x, 5,0);
	 \draw[color =black!100] (-5,\x,0) -- (5, \x,0);
	 }
	 \foreach \x in {-5,...,5}{
	 	\foreach \y in {-5,...,5}{
	\draw (\x,\y,0) [color =black!100] node[scale=.5]{$\bullet$};
		}
         }   
         \draw[color = blue, dashed] (4, 0, 0) -- ( 0,  4, 0)-- ( - 4, 0, 0) -- ( 0,- 4, 0)-- cycle;
         \draw (0,0,0) [color =black!100] node[scale=.7]{$\bullet$};
         \node at (0,-7) {$\Gc_\Mcc\simeq \Z^2$};
         \node[scale=.8] at (0.5,0.5) {$x_\Mcc$};
         	\node[color=blue] at (2,6.5) {$B_{\Gc_\Mcc}(x_\Mcc, r)$};
	\draw[->, color=blue] (3.5,6) .. controls +(-90:.8) and +(45:.2) .. (2.5,2.5);
\end{tikzpicture}
\caption{Localisations at infinity for the $\Z^2$-like graph given in Figure \ref{fig_Z2like}}.
  \label{fig_Z2like2}
\end{figure}
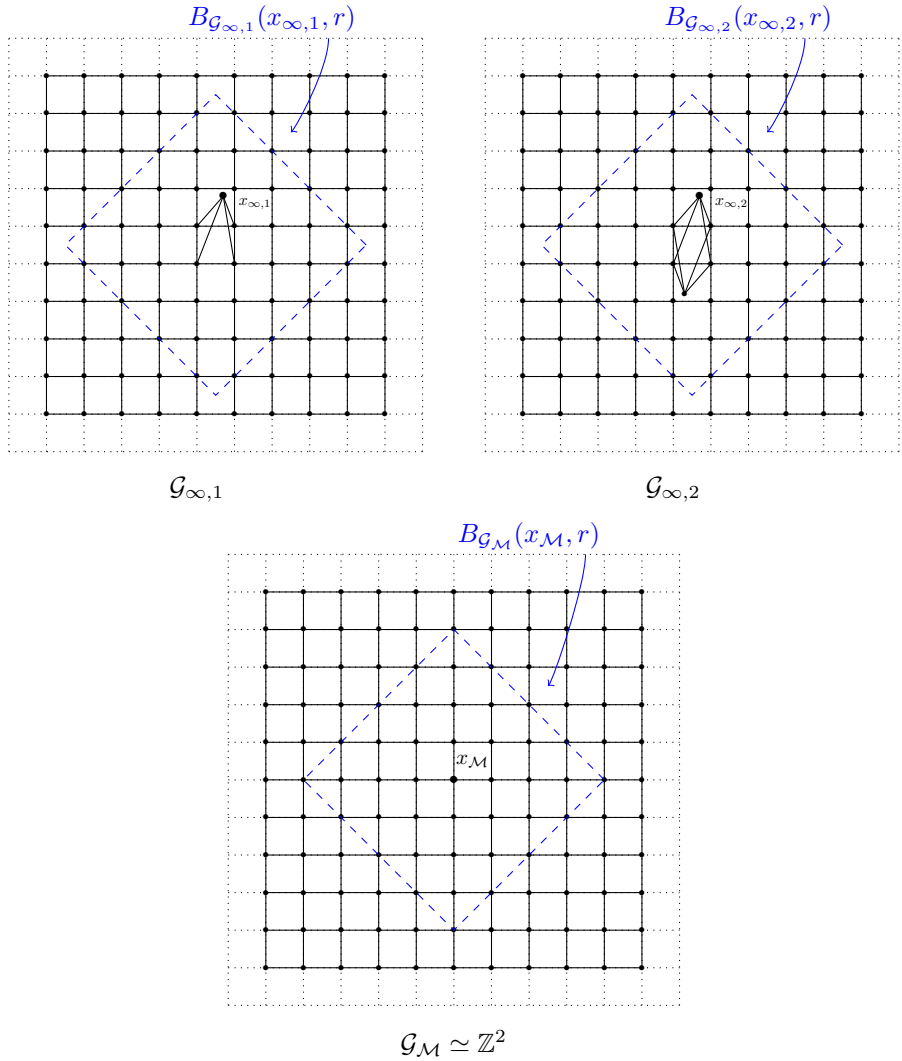

\fi

%%%%%%%%%%%%%%%%%%%%%%%%%%%%%%%%

\section{Proof of the main theorem} 
%In this section we shall prove Theorem \ref{t:simpleF}. 
\subsection{Essential self-adjointness}\label{s:esssa}
In this section we prove the first part of Theorem \ref{t:simpleF}. We rely on a perturbative approach. We assume that $\Gc:= (\Vc, \Ec, m)$ is a Klaus-sparse graph and we use notation of Definition \ref{d:KS2}. We set:
\[\Vc^\sharp:= \bigcup_{i,k} B_\Gc(x_{i,k}, r^{\rm ext}_{i,k}), \quad \Ec^\sharp:= \Ec \times \bone_{\Vc^\sharp\times \Vc^\sharp}, \text{ and } \Gc^\sharp:= (\Vc^\sharp, \Ec^\sharp, m|_{\Vc^\sharp}).\] 
%\brown{Note that $\Gc^\sharp$ is induced by $\Gc$ and that $\|f\|_\Gc= \|f\|_{\Gc^\sharp}$ for all $f$ with support in $\Vc^\sharp$.  }

Due to Definition \ref{d:KS2} (d), the balls that are constituting  $\Vc^\sharp$ are two-by-two disjoint. Using the notation of induced graph given in Section \ref{s:FN}, we  deduce that
\[ \Delta_{\Gc^\sharp} = \bigoplus_{i,k} \Delta_{[B_\Gc(x_{i,k}, r^{\rm ext}_{i,k})]^\Gc}.\]
Since $[B_\Gc(x_{i,k}, r^{\rm ext}_{i,k})]^\Gc$ 
 is a finite graph for all $i,k$, we infer that $\Delta_{\Gc^\sharp}$ is essentially self-adjoint on $\Cc_c(\Vc^\sharp)$. We extend it to $\Vc$ by $0$, we obtain that 
 $\Delta_{\Gc^\sharp} \oplus 0$ is  essentially self-adjoint on $\Cc_c(\Vc)$.

Then, mimicking the proof of \eqref{e:majo}, we have:
\begin{align}\label{e:esssadiff}
|\langle f, (\Delta_\Gc - \Delta_{\Gc^\sharp}\oplus 0) f \rangle_{\Gc}|
\leq 2\langle f, W(\cdot) f \rangle_{\Gc},
\end{align}
where 
\[W(x):= \frac{1}{m(x)}\sum_{y\in \Vc} |\Ec(x,y)- \Ec^\sharp(x,y)|, \text{ for } x\in \Vc.\]
Note that 
\[W(x)= \left\{\begin{array}{ll}
\displaystyle \frac{1}{m(x)}\sum_{y\notin \Vc^\sharp} \Ec(x,y),& \text{ when } x\in \Vc^\sharp,
\\
&
\\
\displaystyle
\frac{1}{m(x)}\sum_{y\in \Vc} \Ec(x,y),& \text{ when } x\notin \Vc^\sharp.
\end{array}\right.\]
Therefore, recalling $(c)$ of Definition \ref{d:KS2}, we have that 
the support of $W$ is contained in $\{\Cc_l\}_{l\in \Lc}$. Moreover, using again (c), we obtain that
\[\sup_{x\in \Vc} W(x)\leq \sup_{x\in \Vc_\Mcc} \deg_{\Gc_\Mcc}(x).\]
In particular,  this implies that $\Delta_\Gc - \Delta_{\Gc^\sharp}\oplus 0$ is a bounded operator. Hence by the Kato-Rellich Theorem, e.g., \cite[Proposition X.12]{RS}, we conclude that $\Delta_\Gc$ is also essentially self-adjoint on $\Cc_c(\Vc)$. \qed

\subsection{General facts about the essential spectrum}\label{s:GFES}
%\subsection{About Weyl sequences}
For the proof of the  second point of the main theorem, we rely extensively on the properties of  approximate eigenfunctions. For the convenience of the reader, we recall these results. 

We start by characterising the spectrum and the essential spectrum of a general self-adjoint operator, e.g., \cite[p.\ 268]{RS}.
\begin{theorem}
Let $H$ be a self-adjoint operator acting in a Hilbert space $(\Hc, \|\cdot\|)$. We have:
\begin{enumerate}[1)]
\item $\lambda \in \sigma(H)$ if and only if there are $\varphi_n \in \Hc$ such that $\|\varphi_n\|=1$
and $\lim_{n\to \infty} (H-\lambda) \varphi_n=0$. The functions $(\varphi_n)_{n\in \N}$ are called \emph{approximate eigenfunctions}.
\item $\lambda \in \sigma_{\rm ess}(H)$ if and only if there are $\varphi_n \in \Hc$ such that 
\begin{enumerate}
\item $\|\varphi_n\|=1$, 
\item $\wlim_{n\to \infty} \varphi_n=0$,
\item $\lim_{n\to \infty} (H-\lambda) \varphi_n=0$. %$\rightharpoonup$
\end{enumerate}
The sequence of functions $(\varphi_n)_{n\in \N}$ is called a \emph{Weyl sequence}.
\end{enumerate}
\end{theorem}

For the next theorem, the first point is a one-line proof. The second one relies on  Persson's Lemma, e.g., \cite[Proposition 18]{KL2}.
\begin{theorem}\label{t:weyl2}
Let $\Gc:=(\Vc,\Ec, m)$ be such that $\Delta_\Gc$ is essentially self-adjoint on $\Cc_c(\Vc)$. It follows that: 
\begin{enumerate}[1)]
\item $\lambda \in \sigma(\Delta_\Gc)$ if and only if there are $\varphi_n \in \Cc_c(\Vc)$ such that $\|\varphi_n\|=1$
and $\lim_{n\to \infty} (\Delta_\Gc-\lambda) \varphi_n=0$.
\item $\lambda \in \sigma_{\rm ess}(H)$ if and only if there are $\varphi_n \in \Cc_c(\Vc)$ such that
\begin{enumerate}
\item $\|\varphi_n\|=1$, 
\item for all finite set $K\subset \Vc$, we have $\supp(\varphi_n)\cap K =\emptyset$, for $n$ large enough,
\item $\lim_{n\to \infty} (\Delta_\Gc-\lambda) \varphi_n=0$. %$\rightharpoonup$
\end{enumerate}
\end{enumerate}
\end{theorem}

\subsection{Computation of the essential spectrum}
In this last section, we finish the proof of the main theorem. We start with two lemmas.

\begin{lemma}\label{l:chicomp}
Let $\Gc:=(\Vc,\Ec, m)$ be such that $\Delta_\Gc$ is essentially self-adjoint on $\Cc_c(\Vc)$. Let $\chi:\Vc\to \R$ be bounded, and define 
\begin{align*}
G(x):= \frac{1}{m(x)}\sum_{y\in \Vc} \Ec(x,y) |\chi(y)-\chi(x)|.
\end{align*}
We have:
\begin{enumerate}
\item For all $f\in \Cc_c(\Vc)$,  we have
\[|\langle f, [\Delta_\Gc, \chi(\cdot)] f\rangle_{\Gc)}| \leq \langle f, G(\cdot) f\rangle_{\Gc}.\]
\item Assuming that $G$ is bounded, $[\Delta_\Gc, \chi(\cdot)]$ extends to a bounded operator that we denote by $[\Delta_\Gc, \chi(\cdot)]_\circ$. 
\item The operator $[\Delta_\Gc, \chi(\cdot)]_\circ$ is compact if 
\begin{align}\label{e:chicomp}
\lim_{|x|\to \infty} G(x)=0.
\end{align}
\end{enumerate}
\end{lemma}
\proof
Take $f\in \Cc_c(\Vc)$. Using the Cauchy-Schwartz inequality, we have:
\begin{align*}
|\langle f, [\Delta, \chi(\cdot)] f\rangle_{\Gc)}| &\leq \sum_{x\in \Vc} |f(x)| \cdot \sum_{y\sim x}  \Ec(x,y) \left|(\chi(y)-\chi(x))\right| \cdot | f(y)|
\\
&\leq 
\frac{1}{2} \sum_{x\in \Vc}\sum_{y\sim x} \Ec(x,y) |\chi(y)- \chi(x)|\cdot |f(x)|^2 + \Ec(x,y)|\chi(y)- \chi(x)|\cdot |f(y)|^2
\\
&=
\sum_{x\in \Vc} G(x) m(x)  |f(x)|^2 =
\langle f, G(\cdot) f\rangle_{\Gc}.
 \end{align*}
 The boundedness is immediate and 
the compactness follows from the min-max theory, e.g., \cite{Gol3}[Propo\-si\-tion 2.8]. \qed

\begin{lemma}\label{l:chiFi}
Assume that  $\Gc:=(\Vc,\Ec, m)$ is Klaus-sparse. Using the notation of Definition \ref{d:KS2}, there exists $\chi:\Vc\to [0,1]$  such that
\begin{enumerate}[1)]
\item $\chi(x)=1$ when $x\in \cup_{i\in \N, k\in \Jc} B_\Gc(x_{i,k}, r_{i,k}^{\rm int})$ and $\chi(x)=0$
when $x\in \left(\cup_{i, k} B_\Gc (x_{i, k}, r_{i,k}^{\rm ext}-2)\right)^c$. 
\item The operator  $[\Delta_\Gc, \chi(\cdot)]_\circ$ is compact.
\end{enumerate}
\end{lemma}
\proof Set
\[\chi(x):= \left\{\begin{array}{rl}
\displaystyle 1- \frac{d_\Gc(x, B_\Gc(x_{i,k}, r_{i,k}^{\rm int}))}{r^{\rm ext}_{i,k} - r^{\rm int}_{i,k}-1},& \text{ if } x\in B_\Gc(x_{i,k}, r_{i,k}^{\rm ext}-1) \text{ for some } i\in \N \text{ and } k\in \Jc,
\\
0,& \text{ otherwise.}
\end{array}\right.\]
Recalling Definition \ref{d:KS2} (d), this can be rewritten as follows:
\[\chi(x)= \sum_{i\in \N, k\in \Jc}\max \left(1- \frac{d_\Gc(x,  B(x_{i,k}, r_{i,k}^{\rm int}))}{r^{\rm ext}_{i,k} - r^{\rm int}_{i,k}-1}, 0 \right) \]
Given $x\in \Vc$, we have
\[G(x):=\frac{1}{m(x)}\sum_{y\sim x} \Ec(x,y) |\chi(y)-\chi(x)| \leq \deg_\Gc(x) \times F(x),\]
where
\[F(x):= \left\{\begin{array}{rl}
\displaystyle \frac{1}{r^{\rm ext}_{i,k} - r^{\rm int}_{i,k}-1},& \text{ if } x\in \cup_{ i\in \N, k\in \Jc} B_\Gc(x_{i,k}, r_{i,k}^{\rm ext}-1)\setminus  B_\Gc(x_{i,k}, r_{i,k}^{\rm int}-1),
\\
0,& \text{ otherwise.}
\end{array}\right.\]
Thanks to the support of $F$ and recalling Definition \ref{d:KS2} (c), we refine the estimate as follows:
\[G(x) \leq \sup_{y\in\Vc_\Mcc} \deg_{\Gc_\Mcc}(y) \times F(x).\]
Moreover, we have that $F(x)\to 0$ as $|x|\to \infty$ due to Definition \ref{d:KS2} (a) and the Fr\'echet convergence. 
This implies \eqref{e:chicomp}, which ensures that
 $[\Delta, \chi(\cdot)]_\circ$ is a compact operator. \qed

\begin{lemma}\label{l:essinfF} Assume that  $\Gc:=(\Vc,\Ec, m)$ is Klaus-sparse. Using the notation of Definition \ref{d:KS2}, 
we have $\sigma(\Delta_\Mcc) \subset \sigma_{\rm ess} (\oplus_{k\in\Jc}\Delta_{\Gc_{\infty, k}})$.
\end{lemma}
\proof Let $\lambda \in \sigma(\Delta_\Mcc)$. First note that $\Delta_\Mcc$ is bounded by \eqref{e:majo} and \eqref{e:Msup}. In particular, $\Cc_c(\Vc_\Mcc)$ is a core for $\Delta_\Mcc$.
Then, by Theorem \ref{t:weyl2}, there are $\varphi_n \in \Cc_c(\Vc_\Mcc)$ so that $\|\varphi_n\|_{\Gc_\Mcc}=1$
and $\lim_{n\to \infty} (\Delta_\Mcc-\lambda) \varphi_n=0$.  Let $r_n>0$ be chosen such that 
$\supp(\varphi_n)\subset B_{\Gc_\Mcc} (x_\Mcc, r_n-1)$. 
 Using properties (d) and (e) of Definition \ref{d:KS2}, we see that there  exist an isometric graph embedding $T_n$  and an injective function $\phi:\N \to \N\times \Jc$ such that 
\[T_n: [B_{\Gc_\Mcc}(x_\Mcc, r_n)]^{\Gc_\Mcc} \to [B_\Gc(x_{\phi(n)}, r^{\rm ext}_{\phi(n)}-1)\setminus  B_\Gc(x_{\phi(n)}, r^{\rm int}_{\phi(n)}+1)]^\Gc \build{\subset}_{}^{(b)} \Gc_{\infty,k(n)},\]
for all $n\in \N$ and $\phi(n)= (i(n), k(n))$. Since $\phi$ is injective, note that $\phi(n)\to \infty$ in the Fr\'echet sense, as $n$ goes to infinity. The last inclusion is due to Definition \ref{d:KS2} (b). 

Since it is a graph embedding and recalling that $\supp(\varphi_n)\subset B_{\Gc_\Mcc} (x_\Mcc, r_n-1)$, we have
\[\left(\oplus_{k\in \Jc}\left(\Delta_{\Gc_{\infty,k}}-\lambda\right)\right) T_n\varphi_n=
(\Delta_{\Gc_{\infty, k(n)}}-\lambda) T_n\varphi_n\build{=}_{}^{(b)}(\Delta_{\Gc}-\lambda) T_n\varphi_n= T_n ( (\Delta_\Mcc-\lambda) \varphi_n) \to 0,\] 
as $n\to\infty$. Moreover, since it is an isometry we have $\|T_n \varphi_n\|_{\cup_k \Gc_{\infty, k}}=1$. Since $\lim_{(i,k)\to \infty}r_{i,k}^{\rm int} = \infty$, up to a subsequence, recalling the support of $T_n\varphi_n$, we see that  $\wlim_{n\to \infty} T_n \varphi_n=0$. This is a Weyl sequence for $\oplus_{k\in \Jc}\Delta_{\Gc_\infty,k}$. This yields that 
$\lambda \in \sigma_{\rm ess}(\oplus_{k\in \Jc}\Delta_{\Gc_\infty,k})$. \qed
\\

We turn to the proof of the main result.
\proof[Proof of Theorem \ref{t:simpleF}:] We prove the two inclusions.

\noindent $\supset $: Let $\lambda \in  \sigma(\Delta_{\Gc_{\infty,k}})$ for some $k\in \Jc$.  One can find functions $\varphi_n \in \Cc_c(\Vc_{\Gc_{\infty,k}})$ such that $\|\varphi_n\|_{\Gc_{\infty, k}}=1$
and $\lim_{n\to \infty} (\Delta_{\Gc_{\infty, k}}-\lambda) \varphi_n=0$ by Theorem \ref{t:weyl2}. Let $r_n>0$ be chosen so that 
$\supp(\varphi_n)\subset B_{\Gc_{\infty, k}} (x_{\infty,k}, r_n-1)$. 
 By 
 %using  $r^{\rm ext}_{i,k} -r^{\rm int}_{i,k} \to \infty$ as $i$ goes to infinity and
  Definition \ref{d:KS2} (b) and (d), we see that there  exist an isometric graph embedding $T_n$ and a strictly increasing function $\phi:\N\to \N$ such that
\[T_n: [B_{\Gc_{\infty, k}}(x_{\infty, k}, r_n)]^{\Gc_{\infty, k}} \to [B_\Gc(x_{\phi(n),k}, r^{\rm int}_{\phi(n),k}-1)]^{\Gc},\]
for all $n\in \N$.  Since it is a graph embedding, we have
\[(\Delta_\Gc-\lambda) T_n\varphi_n=T_n( (\Delta_{\Gc_{\infty, k}}-\lambda) \varphi_n) \to 0,\]
as $n\to \infty$  
and since it is an isometry we have $\|T_n \varphi_n\|_\Gc=1$.

Recalling (d), we see that the supports of $(T_n \varphi_n)_{n\in \N}$ are two by two disjoint. In particular, 
%considering the support of $T_n \varphi_n$, because of Definition \ref{d:KS2} (a), 
we obtain $\wlim_{n\to \infty} T_n\varphi_n=0$. We infer that $T_n\varphi_n$ is a Weyl sequence for $(\Delta_\Gc, \lambda)$. In particular $\lambda \in \sigma_{\rm ess}(\Delta_\Gc)$. This implies that $\cup_{k\in \Jc} \sigma(\Delta_{\infty, k}) \subset \sigma_{\rm ess}(\Delta_\Gc)$. Since $\sigma_{\rm ess}(\Delta_\Gc)$ is closed, we obtain the first inclusion.

\noindent $\subset $: Let $\lambda \in \sigma_{\rm ess}(\Delta_\Gc)$. There are $\varphi_n \in \Cc_c(\Vc)$ verifying 2 (a)-(c) in Theorem \ref{t:weyl2}.  
 Take $\chi$ as in Lemma \ref{l:chiFi} and distinguish two cases.

i) Suppose that $\liminf_{n\to \infty} \|(1-\chi(\cdot))\varphi_n\|_\Gc>0$. Set 
\[\Psi_n:= \frac{1}{ \|(1-\chi(\cdot))\varphi_n\|_\Gc}(1-\chi(\cdot))\varphi_n.\]
We have $\|\Psi_n\|_\Gc=1$, for all $n\in \N$ with support in $ \cap_{i, k} B_\Gc(x_{i,k}, r_{i,k}^{\rm int})^c$. By Definition \ref{d:KS2} (c), we have the following direct sums:
\begin{align}\label{e:suppdel}
\Psi_n = \sum_{l\in \Lc} \bone_{\Cc_l} \Psi_n \text{ and } \Delta_\Gc\Psi_n = \sum_{l\in \Lc} \bone_{\Cc_l} \Delta_\Gc\Psi_n.
\end{align}
For each $l\in \Lc$, we inject $ \bone_{\Cc_l} \Psi_n$ into $\Gc_\Mcc$ and denote it by $\Psi_{n,\Mcc}^l$. Since $\Psi_{n,\Mcc}^{l}$ is with finite support, there is $p=p(l,n)$ such that 
\[1=\sum_{l\in \Lc} \|\Psi_{n,\Mcc}^{l}\|^2_{\Gc_{\Mcc}}= \sum_{l=0}^p \|\Psi_{n,\Mcc}^{l}\|^2_{\Gc_{\Mcc}}.\]
By Lemma \ref{l:amar} (b), there exist $(\theta_{l,n})_{l\in \Lc, n\in \N}$ such that 
\begin{align}\label{e:amarM}
 \left\|\sum_{l\in \Lc} e^{\rmi \theta_{l,n}}\Psi_{n,\Mcc}^l\right\|_{\Gc_{\Mcc}}\geq 1.
\end{align}
Set now
\[ \tilde\Psi_{n, \Mcc}:= \frac{1}{\left\| \sum_{l\in \Lc} e^{\rmi \theta_{l,n}}\Psi_{n,\Mcc}^l  \right\|_{\Gc_{\Mcc}}} \sum_{l\in \Lc} e^{\rmi \theta_{l,n}}\Psi_{n,\Mcc}^l.\]
Then,
%Moreover, up to the embedding given in \eqref{e:backside}, for $r$ large enough, we have:
\begin{align*}
\|(\Delta_\Mcc - \lambda)\tilde \Psi_{n, \Mcc}\|_{\Gc_{\Mcc}}& \build{=}_{}^{}  \frac{1}{\left\| \sum_{l\in \Lc} e^{\rmi \theta_{l,n}}\Psi_{n,\Mcc}^l  \right\|_{\Gc_{\Mcc}}}
\left\|(\Delta_\Gc - \lambda)\sum_{l\in \Lc} e^{\rmi \theta_{l,n}}\bone_{\Cc_l}\Psi_n\right\|_\Gc 
\\
& \build{\leq}_{}^{\eqref{e:amarM}} \left\|\sum_{l\in \Lc} e^{\rmi \theta_{l,n}}(\Delta_\Gc - \lambda)\bone_{\Cc_l}\Psi_n\right\|_\Gc = 
\sum_{l\in \Lc} \left\|(\Delta_\Gc - \lambda)\bone_{\Cc_l}\Psi_n\right\|_\Gc
\\
&\build{=}_{}^{\eqref{e:suppdel}} \left\|(\Delta_\Gc - \lambda)\Psi_n\right\|_\Gc= \frac{\|(\Delta_\Gc - \lambda)(1-\chi(\cdot))\varphi_n\|_\Gc}{  \|(1-\chi(\cdot))\varphi_n\|_\Gc}
\\
&\leq \frac{1}{\|(1-\chi(\cdot))\varphi_n\|_\Gc} \left(\|1-\chi\|_\infty \cdot \|(\Delta_\Gc - \lambda)\varphi_n \|_\Gc + \|[\Delta_\Gc, \chi(\cdot)] \varphi_n\|_\Gc \right) \to 0, 
\end{align*} 
as $n\to \infty$, since $[\Delta_\Gc, \chi(\cdot)]_\circ$ is compact by Lemma \ref{l:chiFi} and $\varphi_n$ tends weakly to $0$. This implies that $\lambda \in \sigma(\Delta_\Mcc)\subset \sigma_{\rm ess} (\Delta_{\oplus_{k\in \Jc}\Gc_{\infty, k}})$, by Lemma \ref{l:essinfF}. 

ii) Suppose now that $\liminf_{n\to \infty} \|(1-\chi(\cdot))\varphi_n\|_\Gc=0$. Up to a subsequence, we can suppose that $\lim_{n\to \infty} \|\chi(\cdot)\varphi_n\|_\Gc=1$. Thanks to Lemma \ref{l:chiFi}, note that
\[\supp (\chi (\cdot)\varphi_n) \subset  \bigcup_{i, k}B_\Gc(x_{i,k}, r_{i,k}^{\rm ext}-2).\]
%Therefore there is $k_0\in \lint 1,j\rint$  such that 
%\[\liminf_{n\to \infty} \|\bone_{\bigcup_{i}B_\Gc(x_{i,k_0}, r_{i,k_0}^{\rm ext})} \chi\varphi_n\| >0\]
Set 
\[\Psi_n :=\frac{1}{\|\chi (\cdot)\varphi_n\|_\Gc}\chi\varphi_n=\frac{1}{\|\bone_{\bigcup_{i,k}B_\Gc(x_{i,k}, r_{i,k}^{\rm ext}-2)} \chi(\cdot)\varphi_n\|_\Gc}\bone_{\bigcup_{i,k}B_\Gc(x_{i,k}, r_{i,k}^{\rm ext}-2)} \chi(\cdot)\varphi_n \in \bigcup_{i,k}B_\Gc(x_{i,k}, r_{i,k}^{\rm ext}-2).\]
Note that $\|\Psi_n\|_\Gc=1$.  Using Definition \ref{d:KS2} (d), we have
\begin{align}\label{e:suppdel2}
\Psi_n= \oplus_{k\in \Jc}\oplus_{i=0}^\infty \bone_{B_\Gc(x_{i, k}, r_{i,k}^{\rm ext}-2)} \Psi_n \text{ and } \Delta_\Gc \Psi_n= \oplus_{k\in \Jc}\oplus_{i=0}^\infty \bone_{B_\Gc(x_{i, k}, r_{i,k}^{\rm ext}-1)} \Delta_\Gc\Psi_n
\end{align}
where the sum is taken over a finite number since $\varphi_n$ is with compact support. Moreover,
\[1=\|\Psi_n\|^2_\Gc= \sum_{k\in \Jc}\sum_{i=0}^\infty \|\bone_{B_\Gc(x_{i, k}, r_{i, k}^{\rm ext}-2)} \Psi_n\|^2_\Gc.\]
We now inject $\bone_{B_\Gc(x_{i, k}, r_{i, k}^{\rm ext}-2)} \Psi_n$ into $\Gc_{\infty, k}$ using Definition \ref{d:KS2} (b). We denote by $\Psi_{n,\infty}^{i,k}: \Vc_{\infty, k} \to \C$ the new function. Trivially, we have
\[\sum_{k\in \Jc}\sum_{i=0}^\infty \|\Psi_{n,\infty}^{i,k}\|^2_{\Gc_{\infty, k}}=1.\]
Since $\Psi_{n,\infty}^{i,k}$ is with finite support, there is $p=p(k,n)$ such that 
\[\alpha_{n,k}:=\sum_{i=0}^\infty \|\Psi_{n,\infty}^{i,k}\|^2_{\Gc_{\infty, k}}= \sum_{i=0}^p \|\Psi_{n,\infty}^{i,k}\|^2_{\Gc_{\infty, k}}.\]
Thanks to Lemma \ref{l:amar} b), there exist $(\theta_{i,k,n})_{i\in \N, k\in \Jc, n\in \N}\in [0,2\pi]^{\N\times \Jc}$ such that 
\begin{align}\label{e:borneinfF}
\sum_{k\in \Jc}  \left\| \sum_{i=0}^\infty e^{\rmi \theta_{i,k,n}} \Psi_{n, \infty}^{i,k}\right\|^2_{\Gc_{\infty, k}}\geq \sum_{k\in \Jc}   \sum_{i=0}^\infty \left\|\Psi_{n, \infty}^{i,k}\right\|^2_{\Gc_{\infty, k}} = \sum_{k\in \Jc} \alpha_{n,k}= 1.
 \end{align}
Set 
\[\tilde \Psi_{n, \infty}^k:= \frac{1}{\sum_{k'\in \Jc}\left\| \sum_{i=0}^\infty e^{\rmi \theta_{i,k',n}} \Psi_{n, \infty}^{i,k'}\right\|_{\Gc_{\infty, k'}}}  \sum_{i=0}^\infty e^{\rmi \theta_{i,k,n}} \Psi_{n, \infty}^{i,k}. \]
Note that $\sum_{k\in \Jc} \|\tilde \Psi_{n, \infty}^k\|_{\Gc_{\infty, k}}=1$ for all $n\in \N$. Using again Definition \ref{d:KS2} (b) and (d), we have
\begin{align*}
\left\|\bigoplus_{k\in \Jc}(\Delta_{\Gc_{{\infty, k}}}- \lambda) \tilde\Psi_{n, \infty}^k\right\|_{\cup_{k\in \Jc}\Gc_{\infty, k}} =& 
\sum_{k\in \Jc} \left\|(\Delta_{\Gc_{{\infty, k}}}- \lambda) \tilde\Psi_{n, \infty}^k\right\|_{\Gc_{\infty, k}}
\\
&\hspace*{-3cm}\build{=}_{}^{(b)} \frac{1}{\sum_{k'\in \Jc}\left\| \sum_{i=0}^\infty e^{\rmi \theta_{i, k',n}} \Psi_{n, \infty}^{i, k'}\right\|_{\Gc_{\infty, k'}}}\sum_{k\in \Jc}\left\|(\Delta_\Gc- \lambda)\left(\sum_{i=0}^\infty e^{\rmi \theta_{i,k,n}}\bone_{B_\Gc(x_{i,k}, r_{i,k}^{\rm ext}-2)} \Psi_n\right) \right\|_\Gc
\\
&\hspace*{-3cm}\build{\leq}_{}^{\eqref{e:borneinfF}}\sum_{k\in \Jc} \left\|\sum_{i=0}^\infty e^{\rmi \theta_{i,k,n}}(\Delta_\Gc- \lambda)\left(\bone_{B_\Gc(x_{i,k}, r_{i,k}^{\rm ext}-2)} \Psi_n\right) \right\|_\Gc 
\\
&\hspace*{-3cm}\build{=}_{}^{(d)}\sum_{k\in \Jc} \sum_{i=0}^\infty \left\| (\Delta_\Gc- \lambda)\left(\bone_{B_\Gc(x_{i,k}, r_{i,k}^{\rm ext}-2)} \Psi_n\right) \right\|_\Gc \build{=}_{}^{\eqref{e:suppdel2}} \|(\Delta_\Gc - \lambda)\Psi_n \|_\Gc
\\
&\hspace*{-3cm}\leq \frac{1}{\|\chi(\cdot)\varphi_n\|_\Gc} \left(\|\chi\|_\infty \cdot \|(\Delta_\Gc - \lambda)\varphi_n \|_\Gc + \|[\Delta_\Gc, \chi(\cdot)] \varphi_n\|_\Gc\right) \to 0, 
\end{align*}
as $n\to \infty$, since $[\Delta_\Gc, \chi(\cdot)]_\circ$ is compact by Lemma \ref{l:chiFi} and $\varphi_n$ tends weakly to $0$. Thus $(\oplus_{k\in \Jc}\tilde\Psi_{n, \infty}^k)_{n\in \N}$ is a Weyl sequence for $\Delta_{\oplus_{k\in \Jc}\Gc_{\infty, k}}$. We conclude that $\lambda \in \sigma(\Delta_{\oplus_{k\in \Jc}\Gc_{\infty, k}})= \overline{\cup_{k\in \Jc} \sigma(\Delta_{\Gc_{\infty, k}})}$. \qed

We have used the following Lemma coming originally from \cite{Am}. 
\begin{lemma}\label{l:amar} a) Let $(e_i)_{i\in \N}$ be included in a Hilbert space $\Hc$ such that $\|e_i\|=1$. Let $(\beta_i)_{i\in \N}\in \C^\N$. %and $(\theta_i)_{i\in \N}\in [0, 2\pi]^\N$. 
For all $p\in \N^*$, we have:
\begin{align}\label{e:amar}
\exists\, \theta_1, \ldots, \theta_p \in [0, 2\pi] \text{ such that } \left\| \sum_{j=1}^p e^{\rmi \theta_j} \beta_j e_j\right\|^2 \geq \sum_{j=1}^p |\beta_j|^2.
\end{align}
b)  In particular, let $(f_j)_{j\in \N}$ be included in a Hilbert space $\Hc$. For all $p\in \N^*$, set  $\alpha(p):=\sum_{j=1}^p\|f_j\|^2$. 
 There exist $(\theta_j)_{j\in \lint 1,p\rint}\subset [0, 2\pi]^{\lint 1,p\rint}$ such that
\[\left\| \sum_{j=1}^p e^{\rmi \theta_j} f_j\right\|^2 \geq \alpha(p).\]
\end{lemma}
\proof a) We prove the result by induction. The case $p=1$ is trivial. Suppose that we have \eqref{e:amar}. Set $k:= \sum_{j=1}^p e^{\rmi \theta_j} \beta_j e_j$. 
Let $\theta_{p+1}\in [0, 2\pi]$ such that $2 \Re \langle e^{\rmi \theta_{p+1}} \beta_{p+1}e_{p+1}, k \rangle \geq 0$. We have:
\begin{align*}
\left\| \sum_{j=1}^{p+1} e^{\rmi \theta_j} \beta_j e_j  \right\|^2 &= \|k\|^2 + |\beta_{p+1}|^2
+ 2 \Re \langle e^{\rmi \theta_{p+1}} \beta_{n+1}e_{p+1}, k \rangle \geq  \|k\|^2 + |\beta_{p+1}|^2.
\end{align*}
which concludes the proof.

b) Set $f_j:= \beta_j e_j$ such that $\|e_j\|=1$ and apply a). \qed
\subsection{Sharpness} We conclude by discussing further examples.
\subsubsection{The closeness of the union is not automatic}\label{s:sharpunion}
In the context of $R$-limits and $C^*$-algebra, the union of the spectra of the localisations at infinity is always closed. A contrario, in the context of Klaus-sparse graphs, we prove in this section that this union is not always a closed set. Theorem \ref{t:simpleF} is sharp. 

To see this let us consider a Klaus-sparse graph $\Gc :=(\Vc,\Ec, m)$ constructed as in Figure \ref{fig_antitree} where the localisations at infinity of $\Gc$ are of the following types:
\begin{enumerate}[1)]
\item Let $\Gc_{\infty, 0}(\Vc_{\infty, 0}, \Ec_{\infty, 0}, 1)$ be a simple $3$-star infinite graph. Namely, let
$\Vc_{{\infty, 0}}:=\{0\} \cup \left(\{1,2,3\}\times \N^*\right)$
and set $\Ec_{{\infty, 0}}(0, (i,1)):=1$ and $\Ec_{\infty, 0}((i, j),
(i, k)):=1$ if $|k-j|=1$, for all $i\in\{1,2,3\}$ and $j,k\in \N^*$ and
set $\Ec_{\infty, 0}(x,y):=0$ otherwise. 
\item
For each $k\in \N^*$, let  $\Gc_{\infty,k} :=( \Vc_{\infty,k},\Ec_{\infty,k}, m_{\infty,k})$ be given by
\[\Vc_{\infty,k} = \Z,  \quad 
\Ec_{\infty,k} (x,y)= \Ec_\Z (x,y) = 
\left\{\begin{array}{ll} 
1,& \text{ if } |x-y|=1,
\\
0, &\text{ elsewhere, }\end{array}\right.
\quad  m_{\infty,k}(x):=
\left\{\begin{array}{ll} 
s(k),& \text{ if } x=0,
\\
1, &\text{ elsewhere.}\end{array}\right.\]
where $s(k)$ is a sequence such that  $0<s(k)<1$ , to be fixed later. 
\item Let $\Gc_\Mcc= (\Z, \Ec_\Z, 1)$.
\end{enumerate}
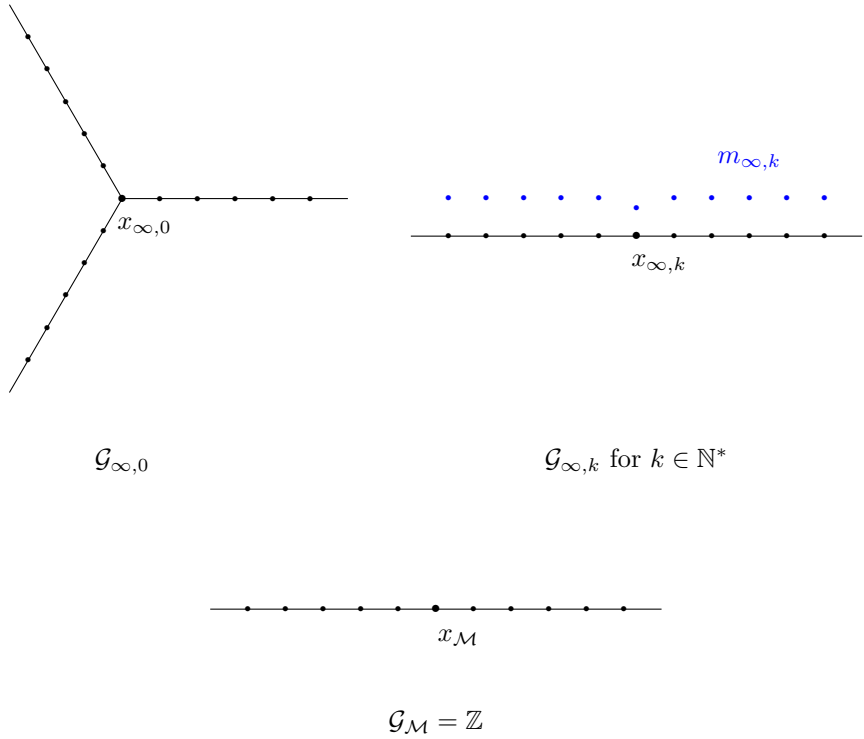
\begin{figure}
%debut G_{\infty,0}
\begin{tikzpicture}[math3d, scale=0.5]
	\draw[color=black](0,0,0)-- (6,0,0);
	\draw[color=black](0,0,0)-- (-6*0.5,6*0.86,0);
	\draw[color=black](0,0,0)-- (-6*0.5,-6*0.86,0);
	 \foreach \x in {0,...,5}{
	 \draw (\x, 0, 0) [color =black!100] node[scale=.5]{$\bullet$};
	 \draw (-\x*0.5, \x*0.86, 0) [color =black!100] node[scale=.5]{$\bullet$};
	 \draw (-\x*0.5, -\x*0.86, 0) [color =black!100] node[scale=.5]{$\bullet$};
	 }
         \node at (0,-7) {$\Gc_{\infty, 0}$};
          \node [color =black!100] at (0.6,-0.75) {$x_{\infty, 0}$};
	 \draw (0, 0, 0) [color =black!100] node[scale=.7]{$\bullet$};
  \end{tikzpicture}
  \quad \quad
% Debut G_{\infty,k}
\begin{tikzpicture}[math3d, scale=0.5]
	%\draw (0,0,0) [color =green] node{$\bullet$};
	\draw[color=black](-6,0,0)-- (6,0,0);
	 \foreach \x in {-5,...,5}{
	 \draw (\x, 0, 0) [color =black!100] node[scale=.5]{$\bullet$};
	 }
	  \foreach \x in {1,...,5}{
	 \draw (\x, 1, 0) [color =blue!100] node[scale=.5]{$\bullet$};
	 \draw (-\x, 1, 0) [color =blue!100] node[scale=.5]{$\bullet$};
	 }
	 \draw (0, 0.75, 0) [color =blue!100] node[scale=.5]{$\bullet$};
	 \node [color =blue!100] at (3,2) {$m_{\infty,k}$};
	 \node [color =black!100] at (0.6,-0.75) {$x_{\infty,k}$};
	 \draw (0, 0, 0) [color =black!100] node[scale=.7]{$\bullet$};
         \node at (0,-6) {$\Gc_{\infty,k}$ for $k\in \N^*$};
\end{tikzpicture}
%debut G_\Mcc
\begin{tikzpicture}[math3d, scale=0.5]
	%\draw (0,0,0) [color =green] node{$\bullet$};
	\draw (0,0,0) [color=white] node[scale=.7]{$\bullet$};
	\draw[color=black](-6,-3,0)-- (6,-3,0);
	 \foreach \x in {-5,...,5}{
	 \draw (\x, -3, 0) [color =black!100] node[scale=.5]{$\bullet$};
	 }
	 \node [color =black!100] at (0.6,-0.75-3) {$x_{\Mcc}$};
	 \draw (0, -3, 0) [color =black!100] node[scale=.7]{$\bullet$};
         \node at (0,-6) {$\Gc_{\Mcc}=\Z$};
\end{tikzpicture}
\caption{Localisations at infinity, counter example for the closeness}
  \label{fig_sharp}
\end{figure}
We refer to Figure \ref{fig_sharp} for an illustration. We have:
\begin{align*}
\sigma(\Delta_{\infty,0})= [0,4]\cup \left\{\frac{9}{2}\right\}, \quad \sigma(\Delta_{\Gc_{\infty, k}})=  [0,4] \cup \left\{\frac{4}{s(k)(2-s(k))}\right\}.
\end{align*}
The former is well-known, e.g, \cite[Lemma 4.4.5]{BrDeEl} with a verbatim proof. The latter comes from a direct computation: the spectrum has an absolute continuous part  and a discrete one, constituted by a unique eigenvalue $\lambda_k=\frac{4}{s(k)(2-s(k))}$.  
Let us choose $s(k):= 1-\frac{\sqrt{2}}{2}+\frac{1}{2k}$. When $k$ goes to $+\infty$,  $(s(k) -1)^2 $ tends to $\frac{1}{2}$, so that $\lambda_k$ tends
to $8$, which shows that the set $\bigcup_{k\in \Jc}\sigma(\Delta_{\Gc_{\infty, k}})$ is not a closed set. A suitable sequence $s_k$ can be chosen in order to get any limit $l$ for $\lambda_k$ , provided $4<l<\infty$.
\subsubsection{A graph without uniform sub-exponential growth}\label{s:nonex} 
Take a graph $\Gc_1:=( \Vc,\Ec, 1))$ such that \eqref{e:growth} is wrong and $\Gc_2:=\Z$. Then take $x\in \Vc$ and glue $\Gc_1$ with $\Gc_2$ in a direct way such that $x$ is identified with $0$. We denote by $\Gc_{\infty,0}$ this graph.  Take a Klaus-graph constructed as in \ref{fig_antitree} where the localisation at infinity is given by $\Gc_{\infty,0}$, see Figure \ref{fig_exp} and $\Gc_\Mcc= \Z$ as medium graph. Clearly $\Gc$ has no uniform sub-exponential growth. The result of \cite{El} cannot apply whereas ours can. 
\begin{figure}
  \def\random{no}              % random sizes of spheres
  \def\f{1.5}                   % strech factor in x-direction
  \def\sizes{1, 2, 6, 24}% (maximal) sizes of spheres
  \newcounter{n}                % n counts the spheres, starting at 0
  \begin{tikzpicture}
      	 \foreach \x in {-9,...,9}{
	 \draw (1.5, \x/4, 0) [color =blue!100] node[scale=.3]{$\bullet$};
	 }
	 \draw [color =blue!100] (1.5, -10/4)--(1.5, 10/4);
    \foreach \j in \sizes
    { \ifthenelse{\equal{\random}{yes}}{
        \pgfmathparse{random(\j)}
        \let\j\pgfmathresult
      }{}
      \ifthenelse{\j>1}{%
        \pgfmathparse{1-\j/2}\let\A\pgfmathresult
        \pgfmathparse{\j/2}\let\B\pgfmathresult
      }{\def\A{0.5}\def\B{1}}%
      \foreach \Y in {\A,...,\B} {
        \coordinate (x) at ({\f*(\then+1)},{(\Y-1/2)/sqrt(\B-\A)}) {};
        \filldraw (x) circle (0.3mm);
        \ifthenelse{\then=0}{}{%
          \foreach \y in {\a,...,\b}
            \draw (\f*\then,{(\y-1/2)/sqrt(\b-\a)}) -- (x);
        }%
      }
      \node [above=1.5mm] [left=0.5mm] at (x) {$S_{\then}$}; % prints S_n
      \stepcounter{n} % next step
      \global\let\a\A % save coordinates for next loop
      \global\let\b\B % save coordinates for next loop
    }
  \end{tikzpicture}
  \addtocounter{n}{-1}
  \ifthenelse{\equal{\random}{yes}}
    {\caption{An antitree with spheres $S_1,\dotsc,S_{\then}$.}}
    {\caption{An antitree with spheres $S_n$ of size $(n+1)!$ glued with $\Z$ at $0$, a counter-example with non-exponential growth}}
  \label{fig_exp}
\end{figure}
\appendix
\section{Stability of the essential spectrum}
Stability of the essential spectrum is a wild subject. The general idea is that if a perturbation is small at infinity then the essential remains the same. To establish it, given $H$ and $H_{\rm pertu}$ being self-adjoint, one usually proves that $(H_{\rm pertu}+\rmi)^{-1} - (H+\rmi)^{-1}$ is compact to obtain that $\sigma_{\rm ess}(H_{\rm pertu})= \sigma_{\rm ess}(H)$. This is the Weyl's theorem, e.g, \cite[Theorem XIII.14]{RS}. The difficulty lies in proving the compactness by taking advantage of the smallness of the perturbation at infinity. We refer  \cite{GeGo2} historical references therein and also for a general abstract method. However, our situation does not fall exactly in the abstract setting of \cite{GeGo2} and we need some adaptations. We present an extended version of \cite[Proposition 5.2]{Gol3} (where $m=\widetilde m$). 

\begin{theorem}\label{t:stabess}
Let $\Gc:=(\Ec, \Vc, m)$ be such that $\Dc\left(\Delta_\Gc^{1/2}\right)= \Dc\left(\deg_\Gc^{1/2}(\cdot)\right)$. Let $\widetilde \Gc:=(\widetilde\Ec, \widetilde\Vc, \widetilde m)$ with $\widetilde \Vc := \Vc$ and such that there is $c$ with $\widetilde \Ec(x,y) \leq c\, \Ec(x,y)$, for all $x,y\in \Vc$ and
\begin{align}\label{e:speed}
\lim_{|x|\to \infty} \frac{m(x)}{\widetilde m (x)}=1 \text{ and } \lim_{|x|, |y|\to \infty} \Ec(x,y) - \widetilde \Ec(x,y) =0.
\end{align}
%where
%\[\Lambda(x):= \frac{1}{m(x)} \sum_{y\in \Vc} |\Ec(x,y)- \widetilde \Ec(x,y)|. \]
Then $\Dc\left(\Delta_{\widetilde \Gc}^{1/2}\right)= \Dc\left(\deg_{\Gc}^{1/2}(\cdot)\right)$ and 
\[\sigma_{\rm ess} (\Delta_{\Gc}) =  \sigma_{\rm ess} (\Delta_{\widetilde \Gc}).\]
\end{theorem}
We stress that the operators $\Delta_{\widetilde \Gc}$ are not necessarily supposed to be essentially self-adjoint $\Cc_c(\Vc)$ and that we consider their Friedrichs's extension. 
\proof 
First we transport unitarily $\Delta_{\widetilde \Gc}$ from $\ell^2(\Vc, \widetilde m)$ into $\ell^2(\Vc, m)$. Namely, we set $\widetilde \Delta:= \sqrt\frac{\widetilde m(\cdot )}{m (\cdot)}\Delta_{\widetilde \Gc} \sqrt{\frac{ m(\cdot )}{\widetilde m (\cdot)}}$. We shall prove that 
\begin{align}\label{e:diffDelta}
\Delta_\Gc- \widetilde \Delta\in \Kc\left( \deg_{ \Gc}^{1/2}(\cdot), \left(\deg_{ \Gc}^{1/2}(\cdot)\right)^*\right),
\end{align} 
here $*$ denotes the antidual and we have identified $\ell^2(\Vc, m)$ with its antidual. 

We have:
\begin{align*}
\Delta_\Gc- \widetilde \Delta&= \left(1- \sqrt\frac{\widetilde m(\cdot )}{m (\cdot)} \right)\Delta_\Gc + \sqrt\frac{\widetilde m(\cdot )}{m (\cdot)} (\Delta_\Gc- \Delta_{\widetilde \Gc})\sqrt{\frac{ m(\cdot )}{\widetilde m (\cdot)}} + \sqrt\frac{\widetilde m(\cdot )}{m (\cdot)} \Delta_\Gc \left(1- \sqrt\frac{ m(\cdot )}{\widetilde m (\cdot)} \right)
\end{align*}
By \eqref{e:majo}, $\Delta_\Gc$ is bounded in $\Bc\left( \deg_{ \Gc}^{1/2}(\cdot), \left(\deg_{ \Gc}^{1/2}(\cdot)\right)^*\right)$ and $ \left(1- \sqrt\frac{\widetilde m(\cdot )}{m (\cdot)} \right) \in \Kc\left( \left(\deg_{ \Gc}^{1/2}(\cdot)\right)^*\right)$ by \eqref{e:speed}. The first term is compact in $\Kc \left( \deg_{ \Gc}^{1/2}(\cdot), \left(\deg_{ \Gc}^{1/2}(\cdot)\right)^*\right)$. In a same way, since $\sqrt\frac{\widetilde m(\cdot )}{m (\cdot)}\in \Bc\left( \left(\deg_{ \Gc}^{1/2}(\cdot)\right)^*\right)$, the third term is also compact. We turn to the second one. Given $f\in\Cc_c(\Vc)$, we have:

\begin{align}
\nonumber
0\leq |\langle f, (\Delta_{\Gc}- \Delta_{\widetilde \Gc}) f\rangle_{\Gc}|&\leq 
\langle f, |\deg_\Gc(\cdot)- \deg_{\widetilde \Gc}(\cdot)| f\rangle_{\Gc}
+\sum_{x,y\in\Vc} |f(x)|\cdot |f(y)|\cdot \left| \Ec(x,y) -\frac{m(x)}{\widetilde m(x)}\widetilde  \Ec(x,y) \right|
\\ \nonumber
&\leq \langle f, |\deg_\Gc(\cdot)- \deg_{\widetilde \Gc}(\cdot)| f\rangle_{\Gc}
+
 \frac{1}{2} \sum_{x,y\in\Vc} |f(x)|^2 \cdot  \left| \Ec(x,y) -\frac{m(x)}{\widetilde m(x)}\widetilde  \Ec(x,y) \right|
 \\ \nonumber
 &\quad \label{e:Fineq}
 + 
  \frac{1}{2} \sum_{x,y\in\Vc} |f(x)|^2\cdot \left| \Ec(x,y) -\frac{m(y)}{\widetilde m(y)}\widetilde  \Ec(x,y) \right|
\\
&\leq \langle f, F(\cdot) f\rangle_{\Gc},
\end{align}
where 
\begin{align}\label{e:F}
F(x)= o(1+\deg_\Gc(x)), \mbox{ as } |x|\to \infty.
\end{align} 
To see this we use \eqref{e:speed} and the Lebesgue convergence. We justify the domination of the first term by
\begin{align*}
\langle f, |\deg_\Gc(\cdot)- \deg_{\widetilde \Gc}(\cdot)| f\rangle_{\Gc} &=
\sum_{x\in\Vc} m(x) |f(x)|^2 \cdot \left| \frac{1}{m(x)}\sum_{y\in\Vc} \Ec(x,y)  
%- \frac{1}{\widetilde m(x)} \sum_{y\in\Vc}  \Ec(x,y) + \frac{1}{\widetilde m(x)} \sum_{y\in\Vc} \widetilde \Ec(x,y) 
- \frac{1}{\widetilde m(x)} \sum_{y\in\Vc} \widetilde \Ec(x,y) \right|
\\
&\hspace*{-2cm}\leq
\sum_{x\in\Vc} m(x) |f(x)|^2 \cdot\left( \left| \left(\frac{1}{m(x)} - \frac{1}{\widetilde m(x)}\right)\sum_{y\in\Vc} \Ec(x,y)\right| + \left| \frac{1}{\widetilde m(x)} \sum_{y\in\Vc} \Ec(x,y)-\widetilde \Ec(x,y) \right|\right)
\\
&
\hspace*{-2cm}\leq
\left(\left(1+\sup_{t\in \Vc} \frac{m(t)}{\widetilde m(t)}\right)+ \left((1+c)\sup_{t\in \Vc} \frac{m(t)}{\widetilde m(t)} \right)\right) \langle f, \deg_\Gc(\cdot) f\rangle
\end{align*}
%%%
and that of the last term, as follows:
\begin{align*}
 \frac{1}{2} \sum_{x,y\in\Vc} |f(x)|^2 \cdot \left| \Ec(x,y) -\frac{m(y)}{\widetilde m(y)}\widetilde  \Ec(x,y) \right|&=   \frac{1}{2} \sum_{x,y\in\Vc} m(x) |f(x)|^2 \cdot \left| \frac{1}{m(x)}\Ec(x,y) - \frac{1}{m(x)}\frac{m(y)}{\widetilde m(y)}\widetilde  \Ec(x,y) \right|
 \\
 &\leq \frac{1}{2}\left(1+ c \sup_{t\in \Vc} \frac{m(t)}{\widetilde m(t)}\right) \langle f, \deg_\Gc(\cdot) f\rangle.
\end{align*}
The treatment of the second term is similar. Thanks to \eqref{e:F}, we infer that the operator $F(\cdot)(1+\deg(\cdot))$ is compact in $\Kc \left( \deg_{ \Gc}^{1/2}(\cdot), \left(\deg_{ \Gc}^{1/2}(\cdot)\right)^*\right)$. By the min-max principle, e.g., \cite[Proposition 2.4]{Gol3}, \eqref{e:Fineq} ensures that the second term also belongs to $\Kc \left( \deg_{ \Gc}^{1/2}(\cdot), \left(\deg_{ \Gc}^{1/2}(\cdot)\right)^*\right)$. We conclude that \eqref{e:diffDelta} follows.
 
We turn to the consequences of \eqref{e:diffDelta}. Using the KLMN's theorem, e.g. \cite[Theorem X.17]{RS}, we obtain that $\Dc(\widetilde\Delta^{1/2})= \Dc(\deg_\Gc^{1/2}(\cdot))$. Going back by unitary transform into $\ell^2(\Vc, \widetilde m)$, we obtain the result for $\Dc(\Delta_{\widetilde \Gc}^{1/2})$. Concerning the equality of the essential spectra, the compactness of $\Delta_\Gc-\tilde \Delta$ implies that $(\Delta_\Gc+\rmi)^{-1}- (\tilde \Delta + \rmi)^{-1}$ is a compact operator, e.g., \cite[Condition (AB)]{GeGo2} or \cite[Proof of Proposition 5.2]{Gol3}. The Weyl's Theorem concludes, e.g, \cite[Theorem XIII.14]{RS}. \qed

To apply this theorem we suppose crucially that the form-domain of $\Delta_\Gc$ is equal to that of $\deg_\Gc(\cdot)$. Recalling \eqref{e:majo}, we have 
$\Dc\left(\deg_\Gc^{1/2}(\cdot)\right)\subset \Dc\left(\Delta_\Gc^{1/2}\right)$ in general but the reverse inclusion is not automatic. We refer to \cite{Gol3} for the beginning of this question and \cite{BoGoKe} for an equivalence. We refer also to \cite{BoGoKeal} for a magnetic version. In our context, it is enough to suppose the equality of the form-domains for the localisations at infinity. This is the aim of the next Proposition.

\begin{proposition}\label{p:fdom}
Let $\Gc:=(\Vc, \Ec, m)$ be a Klaus-sparse graph defined as in Definition \ref{d:KS2}.
Assume that $\Dc\left(\Delta_{\oplus_{k\in \Jc}\Gc_{\infty, k}}^{1/2}\right)=\Dc\left(\deg_{\oplus_{k\in \Jc}\Gc_{\infty, k}}^{1/2}(\cdot)\right)$, for all $k\in \Jc$. Then $\Dc\left(\Delta_{\Gc}^{1/2}\right)=\Dc\left(\deg_{\Gc}^{1/2}(\cdot)\right)$.
\end{proposition}
\proof 
Recalling \eqref{e:majo} it is enough to show that $\Dc\left(\Delta_\Gc^{1/2}\right)\subset\Dc\left(\deg_\Gc^{1/2}(\cdot)\right)$.
By hypothesis and with the help of the uniform boundedness principle, there is $C>0$ such that
\[\langle f, \deg_{\oplus_{k\in \Jc}\Gc_{\infty, k}}(\cdot) f \rangle_{\oplus_{k\in \Jc}\Gc_{\infty, k}} \leq C(\langle f, \Delta_{\oplus_{k\in \Jc}\Gc_{\infty, k}} f \rangle_{\oplus_{k\in \Jc}\Gc_{\infty, k}} +\|f\|^2_{\oplus_{k\in \Jc}\Gc_{\infty, k}}), \]
for all  $f\in \Cc_c(\oplus_{k\in \Jc}\Vc_{\infty, k})$. 
%As in Section \ref{s:esssa} denoting
%\[C_{\Mcc}:=\sup_{x\in \Vc_{\Mcc}} \deg_{\Gc_{\Mcc}}(x)\] 
%such that
%\[\langle f, \deg_{\Gc}(\cdot) f \rangle_{\Gc} \leq (C+C_{\Mcc})(\langle f, \Delta_{\Gc^\sharp}\oplus 0\, f \rangle_{\Gc} +\|f\|^2_{\Gc}), \quad \text{ for all } f\in \Cc_c(\Vc).\]
%Finally using \eqref{e:esssadiff} we obtain there is $c>0$ such that
%\[\langle f, \deg_\Gc(\cdot) f \rangle_\Gc \leq (C+3C_\Mcc)(\langle f, \Delta_\Gc f \rangle_\Gc +\|f\|_\Gc^2), \quad \text{ for all } f\in \Cc_c(\Vc),\]
Take $\chi$ as in Lemma \ref{l:chiFi}. For all $f\in \Cc_c(\Vc)$, we have
\begin{align*}
\langle f, \Delta_{\Gc}\, f \rangle_{\Gc}&= \langle \chi f, \Delta_{\Gc}\, \chi f \rangle_{\Gc} +
\langle \chi f, \Delta_{\Gc}\, (1-\chi)f \rangle_{\Gc}
%\\
%&\quad 
+ \langle (1-\chi)f, \Delta_{\Gc}\, \chi f \rangle_{\Gc}+\langle (1-\chi) f, \Delta_{\Gc}\, (1-\chi) f \rangle_{\Gc}
\\
& \geq \frac{1}{C} \langle \chi f, \deg_\Gc(\cdot)\, \chi f \rangle_{\Gc} - \|\chi f\|_\Gc^2 - 
 \|(1-\chi)\Delta_\Gc\|\cdot \left(2 \|f\|_\Gc\cdot \|\chi f\|_\Gc+ \|f\|_\Gc \cdot \|(1-\chi)f\|_\Gc\right). 
\end{align*}
Recalling that $\|\chi\|=\|1-\chi\|=1$ and that $\|(1-\chi)\Delta_\Gc\|$ and
$\|(1-\chi)\deg_\Gc(\cdot)\|$ are finite by \eqref{e:Msup}, 
we infer there is $c$ such that 
\[\langle f, \deg_\Gc(\cdot) f \rangle_\Gc \leq c(\langle f, \Delta_\Gc f \rangle_\Gc +\|f\|_\Gc^2), \quad \text{ for all } f\in \Cc_c(\Vc),\]
which ensures that $\Dc\left(\Delta_\Gc^{1/2}\right)\subset\Dc\left(\deg_\Gc^{1/2}(\cdot)\right)$ and concludes. \qed


\begin{thebibliography}{xxxxxxxxxx}
\bibitem[Am]{Am} E.\ Amar: \emph{Suites d'interpolation dans le spectre d'une alg\`ebre d'op\'erateurs}, th\`ese d'\'etat, Universit\'e Paris XI, Orsay.
%\bibitem[AF]{AF} C.\ Allard and R.\ Froese: \emph{A Mourre estimate for a
%Schr\"odinger operator on a binary tree}, Rev.\ Math.\ Phys.\ 12
%(2000), no.\ 12, 1655--1667.
%\bibitem[ABG]{ABG} W.~O. Amrein, A.~{Boutet de Monvel}, and  V.~Georgescu: \emph{{$C\sb 0$}-groups,  commutator methods and    spectral theory of {$N$}-body {H}amiltonians},  Progress in Mathematics, vol. 135, Birkh{\"a}user Verlag, Basel, 1996.
%\bibitem[Aom]{Aom} K.\ Aomoto: \emph{Selfadjointness and limit pointness
%for adjacency operators on a tree}, J.\ Analyse Math. 53 (1989), 219--232.
%\bibitem[BoGo]{BoGo} M.\ Bonnefont and S.\ Gol\'enia: \emph{Essential spectrum and Weyl asymptotics for discrete Laplacians}, Ann.\ Fac.\ Sci.\ Toulouse Math. (6) 24 (2015), no.\ 3, 563--624.
\bibitem[BoGoKe]{BoGoKe} M.\ Bonnefont, S.\ Gol\'enia, and M.\ Keller: 
\emph{Eigenvalue asymptotics for Schr\"odinger operators on sparse graphs},  Ann. Inst. Fourier (Grenoble) 65 (2015), no. 5, 1969--1998. 
\bibitem[BoGoKeal]{BoGoKeal} M.\ Bonnefont, S.\ Gol\'enia,  M.\ Keller Liu, Shiping, and Florentin Münch: \emph{Magnetic-Sparseness and Schr\"odinger Operators on Graphs}, Ann.\ Henri Poincar\'e (2020).
\bibitem[BrDeEl]{BrDeEl} J.\ Breuer, S.\ Denisov and L.\ Eliaz: \emph{On the essential spectrum of Schr\"odinger operators on trees}, Math.\ Phys.\ Anal.\ Geom.\ 21 (2018), no.\ 4, Art.\ 33.
%\bibitem[Bir]{Bir} M.\v S.\ Birman: \emph{On the spectrum of singular boundary-value problems}. (Russian)  Mat.\ Sb.\ (N.S.)  55 (97)  1961 125--174. 
%\bibitem[BG]{BG} N.\ Boussaid and S.\ Gol\'enia: \emph{Limiting absorption principle for some long range perturbations of Dirac systems at threshold energies}, Comm.\ Math.\ Phys.\  299  (2010), no. 3, 677--708.
\bibitem[CWL]{CWL} S.N.\ Chandler-Wilde, M.\ Lindner: \emph{Limit Operators, Collective Compactness, and the Spectral theory of Infinite Matrices}, Mem.Amer.Math.Soc. 210, 989 (2011).
\bibitem[Chu]{Chu} F.R.K.\ Chung: \emph{Spectral graph theory}
Regional Conference Series in Mathematics.\ 92.\ Providence, RI:
American Mathematical Society (AMS).\ xi, 207 p.\  
\bibitem[CdV]{CdV} Y.\ Colin de Verdi\`ere:
\emph{Spectres de graphes}, Cours Sp\'ecialis\'es, 4. Soci\'et\'e
Math\'ematique de France, Paris, 1998. 
%\bibitem[CTT]{CTT} Y.\ Colin de Verdi\`ere, N.\ Torki-Hamza and F.\ Truc:
%\emph{Essential self-adjointness for combinatorial Schr\"odinger
%  operators II- Metrically non complete graphs}, Mathematical Physics
%Analysis and Geometry 14, 1 (2011) 21--38. 
%\bibitem[CTT2]{CTT2}Y.\ Colin de Verdi\`ere, N.\ Torki-Hamza and F.\
%  Truc: {Essential self-adjointness for combinatorial Schr\"odinger
%    operators III- Magnetic fields}, Ann. de la Facult\'e des sciences de Toulouse: Math\'ematiques 20 (3), 599-611.
%\bibitem[CDS]{CDS} D.\ Cvetkovi\'c, M.\ Doob, and H.\ Sachs:
%\emph{Spectra of graphs. Theory and application}, Second
%edition. VEB Deutscher Verlag der Wissenschaften, Berlin, 1982.\ 368 pp.
\bibitem[DSV]{DSV} G.\ Davidoff, P.\ Sarnak, and A.\ Valette:
\emph{Elementary number theory, group theory, and Ramanujan graphs},
London Mathematical Society Student Texts, 55.\ Cambridge University
Press, Cambridge, 2003. x+144 pp.
%\bibitem[DG]{DerezinskiGerard} J.~Derezi\'nski and C.~G{\'e}rard: \emph{Scattering theory of classical and quantum $n$-particle systems.}, Texts and Monographs in Physics. Berlin: Springer. xii,444 p., 1997.
%\bibitem[Do]{Do} J.\ Dodziuk: \emph{isoperimetric inequality and
 %   transience of certain random walks},  Trans.\ Amer.\ Math.\ Soc.\ 284
%    (1984),  no.\ 2, 787--794.  
%\bibitem[DK]{DK} J.\ Dodziuk and  W.S.\ Kendall: \emph{Combinatorial
%    Laplacians and isoperimetric inequality}, from local times to
%  global geometry, control and physics (Coventry, 1984/85), 68--74, 
%Pitman Res.\ Notes Math.\ Ser., 150, Longman Sci.\ Tech., Harlow, 1986.
%\bibitem[DM]{DM} J.\ Dodziuk and  V.\ Matthai: \emph{Kato's inequality
%    and asymptotic spectral properties for discrete magnetic
%    Laplacians}, The ubiquitous heat kernel, 69--81, Contemp.\ Math., 398,
%Amer.\ Math.\ Soc., Providence, RI, 2006.
\bibitem[El]{El} L.\ Eliaz: \emph{On the Essential Spectrum of Schrödinger Operators on Graphs}, PhD Thesis, arXiv:1909.10059.
%\bibitem[Fuj]{Fuj} K.\ Fujiwara: \emph{Laplacians on rapidly branching
%    trees},  Duke Math.\ Jour., 83 (1996), No. 1, 191--202.
%\bibitem[GG\'e]{GeorgescuGerard} V.\ Georgescu and C.\ G\'erard: \emph{On the virial theorem in quantum mechanics}, Comm.\ Math.\ Phys.\ 208  (1999) p 275-281.
\bibitem[Ge]{Ge} V.\ Georgescu: \emph{On the structure of the essential spectrum of elliptic operators on metric spaces},  J.\ Funct.\ Anal.\ {\bf 260} (2011), no.\ 6, 1734--1765.
\bibitem[GeGo1]{GeGo1} V.\ Georgescu and S.\ Gol\'enia: \emph{Isometries, Fock spaces, and spectral analysis of
Schr\"odinger operators on trees}, J.\ Funct.\ Anal.\ 227 (2005), no.\ 2, 389--429.
\bibitem[GeGo2]{GeGo2} V.~Georgescu and S.~Gol\'enia: \emph{Decay Preserving Operators and stability of the essential spectrum,}, J.\ Operator Theory  59
%(2008),  no.\ 1, 115--155. 
\bibitem[GeIf1]{GeIf1} V.\ Georgescu and A.\ Iftimovici: \emph{Crossed products of $C^*$-algebras and spectral analysis of quantum Hamiltonians}, Commun.\ Math.\ Phys.\ 228 (2002), no.\ 3, 519--560.
\bibitem[GeIf2]{GeIf2} V.\ Georgescu and A.\ Iftimovici: \emph{Localizations at infinity and essential spectrum of quantum Hamiltonians. I. General theory},  Rev.\ Math.\ Phys.\ {\bf 18} (2006), no.\ 4, 417--483. 
%\bibitem[G\L]{GerardLaba} C.\ G\'erard and I.\ \L aba: \emph{Multiparticle quantum scattering in constant magnetic fields}, Mathematical Surveys and Monographs 90.\ Providence, RI: AMS, American Mathematical Society. xiii. 
\bibitem[Gol1]{Gol1} S.\ Gol\'enia: \emph{$C^*$-algebras of anisotropic Schr\"odinger operators on trees}, J.\ Ann.\ Henri Poincar\'e 5 (2004), no.\ 6, 1097--1115.
\bibitem[Gol2]{Gol2} S.\ Gol\'enia: \emph{Unboundedness of adjacency matrices of locally finite graphs}, Lett.\ Math.\ Phys.\ {\bf 93} (2010), no.\ 2, 127--140.
%\bibitem[Gol2]{Gol2} S.\ Gol\'enia: \emph{Unboundedness of adjacency
%matrices of locally finite graphs}, to appear in Letters in mathematical physics.
\bibitem[Gol3]{Gol3}  S.\ Gol\'enia: \emph{Hardy inequality and asymptotic eigenvalue distribution for discrete Laplacians}, 
J.\ Funct.\ Anal.\ {\bf 266} (2014), no.\ 5, 2662--2688.
%\bibitem[GS]{GS} S.\ Gol\'enia and C.\ Schumacher: \emph{The problem
%    of deficiency     indices for discrete Schr\"odinger operators on locally finite graphs}, to appear in Journal of Mathematical Physics. 
%\bibitem[GJ]{GoleniaJecko} S.~Gol{\'e}nia and T.~Jecko: \emph{{A new look at Mourre's commutator  theory.}}, Complex Anal.\ Oper.\ Theory \textbf{1} (2007), no.~3, 399--422.
%\bibitem[GM]{GoleniaMoroianu} S.~Gol{\'e}nia and S.~Moroianu: \emph{Spectral analysis of magnetic Laplacians on conformally cusp manifolds}, Ann. Henri Poincar{\'e} 9 (2008), no.~1, 131--179.
%\bibitem[HS]{HelfferSjostrand} B.~Helffer and J.~Sj{\"o}strand: \emph{Op{\'e}rateurs de schr{\"o}dinger avec champs magn{\'e}tiques faibles et constants}, S{\'e}min. {\'E}quations  D{\'e}riv. Partielles 1988--1989, Exp.\ No.\ 12, 11 p.\ (1989).
%\bibitem[HSS]{HunzikerSigalSoffer} W.\ Hunziker, I.M.\ Sigal, and A.\ Soffer: \emph{Minimal escape velocities}, Communications in Partial Differential Equations, 24:11, 2279 --2295 (1999). 
%\bibitem[Jor]{Jor08} P.E.T.\ J\o rgensen: \emph{Essential
%  self-adjointness of the graph-Laplacian}  J. Math. Phys. 49, No.\ 7, 073510, 33 p.(2008).
%\bibitem[JP]{JP} P.E.T.\ J\o rgensen and E.P.J.\ Pearse:
%  \emph{Spectral reciprocity and matrix representations of unbounded
 %   operators}, to appear in Journal of Functional Analysis. 
\bibitem[Kla]{Kla} M.\ Klaus: {On $-d^{2}/dx^{2}+V$ where
    $V$ has infinitely many ``bumps''},
Ann.\ Inst.\ H.\ Poincar\'e Sect.\ A (N.S.) 38 (1983), no.\ 1, 7--13. 
\bibitem[Kel]{Kel} M.\ Keller: \emph{The essential spectrum of the
    Laplacian on rapidly branching tessellations}, Math.\ Ann.\ 346,
  Issue 1 (2010), 51--66. 
\bibitem[KL]{KL2} M.\ Keller and D.\ Lenz: \emph{Unbounded Laplacians
    on graphs: Basic spectral properties and the heat equation},
  Math.\ Model.\   Nat.\ Phenom.\ (2009) Vol.\ 5, No.\ 2.
\bibitem[KL2]{KL} M.\ Keller and D.\ Lenz: \emph{Dirichlet forms and
  stochastic completeness of graphs and subgraphs}, to appear in
Journal fuer die reine und angewandte Mathematik.
%\bibitem[KLW]{KLW} M.\ Keller, D.\ Lenz and R.\ Wojciechowski:
 % \emph{Volume Growth, Spectrum and Stochastic Completeness of
  %  Infinite Graphs}, preprint arXiv:1105.0395v1. 
\bibitem[LaSi]{LaSi} Y.\ Last and B.\ Simon: \emph{The Essential Spectrum of Schr\"odinger, Jacobi, and CMV Operators}, J.\ Anal.\ Math.\ 98 (2006), 183--220.
\bibitem[M\u{a}n]{Man} M.\ M\u{a}ntoiu: \emph{$C^*$-algebras, dynamical systems at infinity and the essential spectrum of generalized Schr\"odinger operators}, J.\ Reine Angew. Math.\ 550 (2002), 211--229.
\bibitem[MaPuRi]{MaPuRi} M.\ M\u{a}ntoiu, R.\ Purice, and S.\ Richard: \emph{Spectral and propagation results for magnetic Schr\"odinger operators; a $C^*$-algebraic framework}, J.\ Funct.\ Anal.\ 250 (2007), no.\ 1, 42--67.
%\bibitem[Mas]{Mas} J.\ Masamune: \emph{A Liouville property and its
%    application to the Laplacian of an infinite graph}, Spectral
%    analysis in geometry and number theory,  103--115, Contemp.\ Math., 484, Amer.\ Math.\ Soc.
%\bibitem[Mil]{Mil} O.\ Milatovic: {\em Essential self-adjointness of
 %   magnetic Schr\"odinger operators on locally finite graphs}, to appear in Integral Equations and Operator Theory.
%\bibitem[Mil2]{Mil2} O.\ Milatovic: {\em Essential self-adjointness of
%    discrete magnetic Schr\"odinger operators}, preprint arXiv:1105.3129v1.
%\bibitem[MW]{MW} B.\ Mohar and W.\ Woess: {\em A survey on spectra of infinite graphs}, J.\ Bull.\ Lond.\ Math.\ Soc.\ 21, No.3, 209-234 (1989).
\bibitem[NaTa1]{NaTa1} S.A.\ Nazarov and J.\ Taskinen: \emph{Essential spectrum of a periodic waveguide with non-periodic perturbation}, J.\ Math.\ Anal.\ Appl.\ {\bf 463} (2018), no.\ 2, 922--933.
\bibitem[NaTa2]{NaTa2} S.A.\ Nazarov and J.\ Taskinen: \emph{Essential spectrum of periodic medium with sparsely placed foreign inclusions}, preprint. 
\bibitem[Ra]{Ra} V.S.\ Rabinovich: \emph{The essential spectrum of Schr\"odinger operators on periodic graphs} (Russian) Funktsional.\ Anal.\ i Prilozhen.\ {\bf 52} (2018), no.\ 1, 80--84; translation in Funct.\ Anal.\ Appl.\ {\bf 52} (2018), no.\ 1, 66--69.
\bibitem[RS]{RS}M.\ Reed and B.\ Simon: \emph{ Methods of Modern Mathematical
Physics, Tome I--IV: Analysis of operators} Academic Press.
%\bibitem[Sch]{Sch} J.\ Schwinger: \emph{On the bound states of a given
 %   potential}, Proc.\ Nat.\ Acad.\ Sci.\ U.S.A., 47:122--129, 1961.
\bibitem[S]{S}B.\ Simon, \emph{Szego's Theorem and its Descendants}, M. B. Porter Lectures, Princeton
University Press, Princeton, NJ, 2011. Spectral Theory for L2 Perturbations of Orthogonal Polynomials.
\bibitem[SaSu]{SaSu} I.\ Sasaki and A.\ Suzuki: \emph{Essential spectrum of the discrete Laplacian on a perturbed periodic graph},  J.\ Math.\ Anal.\ Appl.\ {\bf 446} (2017), no.\ 2, 1863--1881.
%\bibitem[Tor]{Tor} N.\ Torki-Hamza: \emph{Laplaciens de graphes
%    infinis I Graphes m\`etriquement complets}, Confluentes Math.\ 2 (3) (2010) 333--350.
%\bibitem[Web]{Web} A.\ Weber: \emph{Analysis of the physical Laplacian
 %   and the heat flow on a locally finite graph}, J.\ Math.\ Anal.\  Appl.\ 370 (1) (2010) 146--158.
\bibitem[Woj]{Woj} R.\ Wojciechowski: {\em Stochastic
completeness of graphs}, Ph.D.\ Thesis, 2007, arXiv:0712.1570v2[math.SP].
\end{thebibliography}
\end{document}